\documentclass{article}

\usepackage{a4}
\usepackage{paralist}

\title{On the importance of smoothness, interface resolution and numerical sensitivities in shape and topological sensitivity analysis}
\author{M.H. Gfrerer \\
  \emph{Institute of Applied Mechanics, Graz University of Technology} \\
  P. Gangl \\
  \emph{Johann Radon Institute for Computational and Applied }\\\emph{Mathematics (RICAM), Austrian Academy of Sciences}}

\usepackage{amsthm}
\usepackage{amsfonts}
\usepackage{booktabs}
\usepackage{graphicx}
\usepackage{diagbox}

\usepackage{multirow}
\usepackage{hyperref}
\usepackage{xcolor}

\usepackage[english]{babel}  
\usepackage[intlimits]{amsmath}
\usepackage{amsthm}
\usepackage{amssymb} 
\usepackage{units} 
\usepackage{amsfonts,mathrsfs,bbm}
\allowdisplaybreaks[1]
\usepackage[numbers]{natbib}
\usepackage{graphicx}
\usepackage{tikz}
\usepackage[T1]{fontenc}
\usepackage{mathptmx}
\usepackage[scaled=.92]{helvet}
\usepackage{courier}
\usepackage{tikz-3dplot}
\usetikzlibrary{shapes,arrows}
\usepackage{pgfplots}
\pgfplotsset{every axis/.append style={thick}}
\pgfplotscreateplotcyclelist{mylines}{%
  {red,mark=*},
  {blue,mark=square},
  {green,mark=+},
  {black,mark=o},
  {violet,mark=otimes},
  {brown,mark=triangle},
  {cyan,mark=diamond},
  {orange,mark=10-pointed star},
  {magenta,mark=|}  } 
\usetikzlibrary{external}

\usepackage{standalone}
\usepackage{hyperref}
\usepackage{cleveref}
\usepackage{nicefrac}
\usepackage{subcaption}
\usepackage{epstopdf}
\usepackage{color}
\usepackage{mathtools}
\usepackage{tikz}
\usepackage{pgfplots}
\usepackage{booktabs}
\usepackage{array,makecell}
\usepackage{mdframed}
\usepackage{xspace}
\usetikzlibrary{calc,shapes}

\pgfplotsset{compat=1.13}
\usetikzlibrary{math}
\usetikzlibrary{intersections}
\usetikzlibrary{arrows.meta}
\newtheorem{theorem}{Theorem}[section]

\newtheorem{lemma}[theorem]{Lemma}
\newtheorem{remark}[theorem]{Remark}

\definecolor{col1}{RGB}{0, 114, 254}   
\definecolor{col2}{RGB}{255, 159, 0}  
\definecolor{col3}{RGB}{0, 155, 0}   
\definecolor{col4}{RGB}{213, 94, 0}  
\definecolor{col5}{RGB}{152,78,163}  
\definecolor{col6}{RGB}{230,14,80}  

\usepackage[ulem=normalem,draft]{changes}

\newcommand*\Heq{\ensuremath{\overset{\kern2pt H}{=}}}
\newcommand{\FEM}{\textit{FEM}\xspace}
\newcommand{\XFEM}{\textit{XFEM}\xspace}

\newcommand{\topshapeDerivative}{topological-shape\xspace}

\newcommand{\corrections}[1]{{#1}}

\newcommand{\targetU}{\hat u}
\newcommand{\uu}{\mathbf u}
\newcommand{\uuzero}{{\uu}}
\newcommand{\pp}{\mathbf p}

\newcommand{\KK}{\mathbf K}
\newcommand{\MM}{\mathbf M}
\newcommand{\mm}{\mathbf m}
\newcommand{\Amat}{\mathbf K}

\newcommand{\MmatSEeps}{{\mathbf M}^{SE}_\eps}

\newcommand{\MmatSE}{{\mathbf M}^{SE}_0}

\newcommand{\KmatSeps}{{\mathbf K}^S_\eps}
\newcommand{\KmatSEeps}{{\mathbf K}^{SE}_\eps}
\newcommand{\KmatESeps}{{\mathbf K}^{ES}_\eps}

\newcommand{\KEeps}{{K}^E_\eps}

\newcommand{\KmatSE}{{\mathbf K}^{SE}_0}

\newcommand{\KE}{{K}^E_0}
\newcommand{\Q}{{\mathbf Q}_\eps}
\newcommand{\Qtilde}{{\mathbf Q}^\top_\eps}
\newcommand{\ff}{\mathbf f}

\newcommand{\ie}{\textit{i.e.}\,}

\newcommand{\R}{\mathbb R}

\newcommand{\ddx}{\;\mathrm{d} { x}}
\newcommand{\eps}{\varepsilon}
\newcommand{\objectiveFunction}{G}

\newcommand{\redObjectiveFunction}{{\hat{G}}} 

\newcommand{\enrichmentBasis}{N^{E}}

\newcommand{\standardMethod}{\textit{standard method}\xspace }
\newcommand{\enrichedMethod}{\textit{enriched method}\xspace }
\newcommand{\pertubationTopology}{\eps}

\newcommand{\characteristic}{\chi}

\graphicspath{{pics/}{pic/}{code/results/}}

\newtheorem{problem}{Problem}

\numberwithin{equation}{section} 

\begin{document}

\maketitle

\begin{abstract}
In this paper we investigate the influence of the discretization of PDE constraints on shape and topological derivatives. To this end, we study a tracking-type functional and a two-material Poisson problem in one spatial dimension. We consider the discretization by a \textit{standard method} and an \textit{enriched method}. In the standard method we use splines of degree $p$ such that we can control the smoothness of the basis functions easily, but do not take any interface location into consideration. This includes for $p=1$ the usual hat basis functions. In the enriched method we additionally capture the interface locations in the ansatz space by enrichment functions. For both discretization methods shape and topological sensitivity analysis is performed.
It turns out that the regularity of the shape derivative depends on the regularity of the basis functions. Furthermore, for point-wise convergence of the shape derivative the interface has to be considered in the ansatz space. For the topological derivative we show that only the enriched method converges.
\end{abstract}

%
%
%



%
%
\tableofcontents
\section{Introduction}

There are a number of different methods for shape and topology optimization, which might be roughly grouped into Lagrangian methods, density based methods, and level-set methods.
In Lagrangian methods, shapes are represented using a computational mesh. The advantage of this approach lies in its ability to provide accurate mechanical analysis, as long as the mesh quality remains acceptable. However, a significant drawback is that the mesh quality deteriorates during optimization, necessitating re-meshing, which is notoriously difficult. Additionally, Lagrangian methods cannot accommodate topological changes, meaning the creation or closure of holes is not possible.

In the class of density-based topology optimization methods \cite{bendsoe2003topology}, a design is represented by a density function $\rho(\mathbf x)$ that is allowed to attain any value in the interval $[0,1]$. Then, regions with $\rho(\mathbf x) = 0$ and $\rho(\mathbf x) = 1$ are interpreted as occupied by material 1 and 2, respectively, while intermediate density values $0<\rho(\mathbf x) < 1$ are penalized in order to obtain designs that are almost 'black-and-white'. One advantage of density based methods is that the system response depends continuously on $\rho$ and the standard notions of derivatives in vector spaces can be applied. While topological changes are possible, interfaces are typically not crisp and there is no measure of optimality with respect to shape variations at the interface.

Due to the shortcomings of Lagrangian methods and density-based methods, level-set method have been developed. Here, shapes are represented in an implicit way by a level set function $\phi(\mathbf x)$. This offers sharp designs and the possibility for topological changes in a natural way. However, if the level-set function is driven by shape sensitivities only it lacks a nucleation mechanism \cite{allaire2004structural}. Then one has to use a perforated initial design for topology optimization problems. This drawback can be overcome by using topological sensitivity information. One possibility is to use alternating shape or topological update steps \cite{AllaireJouve:2006a,burger2004incorporating}, while in \cite{amstutz2006new} an algorithm solely based on topological derivatives is presented. However, that latter method uses an average of topological derivatives as sensitivities at the material interfaces, which does not necessarily represent the correct shape sensitivities.
In \cite{gangl2023unified} we have introduced the \topshapeDerivative derivative, which unifies shape and topological derivatives, such that in the subdomains topological sensitivity information and at interfaces shape sensitivity information is combined. This allows simultaneous shape and topological update steps, without averaged sensitivities. In level-set methods one further issue is the accurate mechanical analysis. Typically, simple interpolation procedures within cut elements are used \cite{amstutz2006new,gangl2020multi,gangl2023unified}, which do not properly resolve the interface and lead to "rough" optimal designs. Such non-optimal behavior due to non-resolved interfaces is not exclusive to optimization problems, but a prevalent issue encountered in numerical simulations, like multi-phase flows, plasticity problems, or contact problems. In order to migrate the drawbacks of unresolved interfaces one can (i) use re-meshing \cite{allaire2014shape}, or (ii) use a smooth discretization space, or (iii) use an unfitted finite element method. The latter two approaches are investigated in the present paper. Motivated by numerical experience from computational contact mechanics \cite{temizer2012three,corbett2015three} we investigate the use of smooth basis functions. In the cited works the use of B-splines reduces oscillations compared to classical Lagrange polynomials and thus improves numerical stability. In order to avoid re-meshing also advanced "unfitted finite element methods" like XFEM, enriched FEM and CutFEM with optimal convergence properties have emerged \cite{barrett1987fitted,moes1999finite,BurmanClausHansboLarsonMassing2014,burman2025cut}. We study also the use of an enriched method for the evaluation of sensitivities.

Another longstanding question in optimal design concerns whether sensitivity information should be derived from the continuous problem or from its discretized counterpart. The former corresponds to the \textit{optimize-then-discretize} approach, which is extensively discussed in \cite{DZ2,hinze2008optimization,henrot2018shape,allaire2021shape}, while the latter is known as the \textit{discretize-then-optimize} approach. A unified framework that bridges the continuous and discrete settings is presented in \cite{berggren2009unified}.
The \textit{discretize-then-optimize} approach yields exact sensitivities for the discretized problem \cite{glowinski1998shape}, but it must be derived separately for each finite element formulation. In contrast, the \textit{optimize-then-discretize} approach is independent of the discretization scheme, although it provides correct sensitivity information only for the exact solution.
Finally, we note that shape sensitivity analysis for discretized problems has been developed for XFEM discretizations in \cite{noel2016analytical} and for CutFEM discretizations in \cite{berggren2023shape,ShahanWalker2025,wegert2025level}.

In the present paper we analyze the \textit{discretize-then-optimize} approach for different discretization schemes and shape as well as topological perturbations. In the following we call the sensitivities obtained in this way the numerical shape and topological sensitivities, respectively. In particular, we investigate the influence of
\begin{enumerate}
	\item the smoothness of the discretization,
	\item and the interface resolution of the discretization,
\end{enumerate}
on the properties of the
\begin{enumerate}
\item numerical shape sensitivities,
\item and numerical topological sensitivities.
\end{enumerate}

Moreover, we compare these numerical sensitivities with the continuous sensitivities obtained through the \textit{optimize-then-discretize} approach. To enable a direct comparison with the fully analytical sensitivities, we consider a tracking-type functional associated with a two-material Poisson problem in one spatial dimension. Different discretization variants are examined, employing basis splines of arbitrary order, optionally enhanced by global linear enrichment functions. The use of splines provides convenient control over the smoothness of the discretization.
While the choice of linear enrichment functions yields suboptimal convergence rates for higher order splines, we chose this problem setup in order to allow for explicit computation of the mentioned numerical sensitivities, allowing us to investigate the effect of resolving the interface by the discretization.

This article is organized as follows. The next section presents the model problem and
basic material from shape and topological sensitivity analysis in the continuous setting. Furthermore, analytic solutions and sensitivities are computed explicitly. Section 3  describes the different discretization methods studied in the paper. Then, in Section 4 shape and topological sensitivity analysis is performed for the discretized problem. A theorem regarding the smoothness of the shape derivative is given in Section 5.  Numerical experiments are discussed in Section 6. Eventually, we draw some conclusions around the present study in Section 7 and outline a few possible topics for future work.

\section{Problem settings and continuous shape and topological derivative}
\label{sec::continuous}
In the present paper we restrict ourselves to 1D problems, \ie the hold-all domain is $D = (0,\ell)$ with $\ell>0$ and as a PDE constraint we consider the Poisson equation. We search for a partitioning $\Omega_1$, $\Omega_2$, such that $D = \Omega_1 \cup \Omega_2$ and $\Omega_1 \cap \Omega_2 = \emptyset$ holds. For such a partition, we introduce the symbol $\Omega := (\Omega_1, \Omega_2)$. With $\Omega_1$ and $\Omega_2$ we associate the material parameters $\lambda_1$ and $\lambda_2$ respectively. The precise problem formulation is given next.
\begin{mdframed}
\begin{problem}\label{prob::tracking}
	For fixed parameters $\ell,\lambda_1, \lambda_2>0$, the source term $f(x):D\rightarrow\R$ and the desired state $\targetU:D\rightarrow\R$ the optimization problem reads:
	\begin{subequations}\label{eq::1dproblemSettingTracking}
		\begin{alignat}{2}
			\min_{\Omega_2\in \mathcal A} \objectiveFunction(u) &= \int_{0}^{\ell} \big(u(x) - \targetU(x)\big)^2 \ddx \\
			\text{ subject to}& \nonumber \\
			-\frac{d}{dx}\left(\lambda_\Omega\frac{du(x)}{dx}\right) &= f(x) \quad\text{for } x\in D, \label{eq::trackingA}\\
			u(0) &= 0, \label{eq::trackingB}\\
			u(\ell) &= 0 \label{eq::trackingC},\\
			\left(\lambda_\Omega\frac{du}{dx}\right)\bigg|_{\zeta^-} &= \left(\lambda_\Omega\frac{du}{dx}\right)\bigg|_{\zeta^+}\label{eq::trackingD} \quad \forall \zeta \in \overline{\Omega}_1\cap \overline{\Omega}_2.
		\end{alignat}
	\end{subequations}
\end{problem}
\end{mdframed}
Here, $\mathcal A$ is the set of admissible shapes and the piecewise constant diffusion coefficient $\lambda_\Omega$ is given by
\begin{align*}
	\lambda_\Omega(x) =& \characteristic_ {\Omega_1}(x)\lambda_1 + \characteristic_{\Omega_2}(x) \lambda_2,
\end{align*}
with $\characteristic_S$ the characteristic function of a set $S$,
\begin{align*}
	\characteristic_S(x) = \begin{cases}
		1, & x \in S, \\
		0, & \mbox{ else}.
	\end{cases}
\end{align*}
Let $\Gamma_D = \{0,\ell\}$ and $H^1_0(D) = \{v \in H^1(D): v = 0 \mbox{ on } \Gamma_D \}$. The weak formulation of the PDE constraint \labelcref{eq::trackingA,eq::trackingB,eq::trackingC,eq::trackingD} reads
\begin{align} \label{eq_weakForm}
	&\mbox{Find }u\in H_0^1(D)
    \mbox{ such that } \nonumber \\
	&\int_D \lambda_\Omega \frac{du}{dx} \frac{d v}{dx} \ddx = \int_D f(x)\, v \ddx \quad  \mbox{ for all }v \in H_0^1(D).
\end{align}
For a given domain $\Omega_2\in \mathcal A$, \eqref{eq_weakForm} admits a unique solution, which we denote by $u(\Omega_2)$. This allows to introduce the shape function $\redObjectiveFunction(\Omega_2):= \objectiveFunction( u(\Omega_2))$.
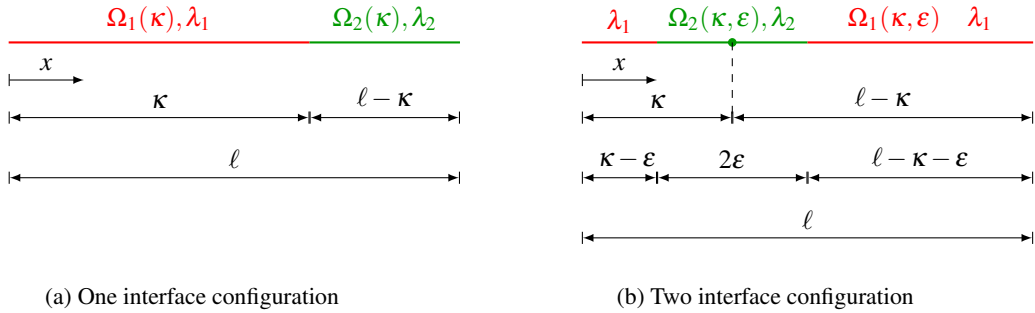
\begin{figure}[ht]
	\begin{subfigure}[t]{0.4\textwidth}
		\begin{tikzpicture}[>=latex]
            \useasboundingbox (-0.1,-3) rectangle (6.1,1);
			\draw[thick,red] (0,0) -- (4,0)node[pos=0.5,above]{$\Omega_1(\kappa), \lambda_1$};
			\draw[thick,green!60!black] (4,0) -- (6,0)node[pos=0.5,above]{$\Omega_2(\kappa), \lambda_2$};
			\draw[|->] (0,-0.5) -- ++(1,0) node[pos=0.5,above]{$x$};
			\draw[|<->|] (0,-1) -- ++(4,0) node[pos=0.5,above]{$\kappa$};
			\draw[|<->|] (4,-1) -- ++(2,0)node[pos=0.5,above]{$\ell-\kappa$};
			\draw[|<->|] (0,-1.8) -- ++(6,0) node[pos=0.5,above]{$\ell$};
		\end{tikzpicture}
		\caption{One-interface configuration \corrections{as used for shape perturbations}}
		\label{fig:flowOneInterface}
	\end{subfigure}
	\hfill
 \begin{subfigure}[t]{0.4\textwidth}
		\begin{tikzpicture}[>=latex]
            \useasboundingbox (-0.1,-3) rectangle (6.1,1);
			\draw[thick,red] (0,0) -- (1,0)node[pos=0.5,above]{$\lambda_1$};
			\draw[thick,green!60!black] (1,0) -- (3,0)node[pos=0.5,above]{$\Omega_2(\kappa,\eps), \lambda_2$};
            \draw[fill,green!60!black] (2,0) circle(0.05);
            \draw[dashed,thin] (2,0) -- ++(0,-1);
			\draw[thick,red] (3,0) -- (6,0)node[pos=0.5,above]{$\Omega_1(\kappa,\eps)\quad\lambda_1$};
			\draw[|->] (0,-0.5) -- ++(1,0) node[pos=0.5,above]{$x$};
            \draw[|<->|] (0,-1) -- ++(2,0) node[pos=0.5,above]{$\kappa$};
			\draw[|<->|] (2,-1) -- ++(4,0)node[pos=0.5,above]{$\ell-\kappa$};
			\draw[|<->|] (0,-1.8) -- ++(1,0) node[pos=0.6,above]{$\kappa-\pertubationTopology$};
			\draw[|<->|] (1,-1.8) -- ++(2,0) node[pos=0.5,above]{$2\pertubationTopology$};
			\draw[|<->|] (3,-1.8) -- ++(3,0) node[pos=0.5,above]{$\ell-\kappa-\pertubationTopology$};
			\draw[|<->|] (0,-2.6) -- ++(6,0) node[pos=0.5,above]{$\ell$};
		\end{tikzpicture}
		\caption{Two-interface configuration \corrections{as used for topological perturbations}}
		\label{fig:flowTwoInterface}
	\end{subfigure}
	\caption{1D problem setups}
	\label{fig:problemillustration1d}
\end{figure}

\subsection{Continuous shape derivative}\label{sec::ContinuousShapeDerivative}
In this section, we introduce a continuous shape sensitivity, which we investigate in the present paper. For simplicity we consider only configurations with one interface at the position $\kappa\in (0,\ell)$, such that $\Omega_1(\kappa) = (0,\kappa)$ and $\Omega_2(\kappa)=(\kappa,\ell)$. Such a situation is illustrated in \Cref{fig:flowOneInterface}.
This has the advantage that we can compute analytic solutions to the state and the adjoint problem, which depend symbolically on the interface location $\kappa$. Thus, it is possible to symbolically differentiate the objective function. Therefore, the shape derivative for the considered configuration is defined by
\begin{align}\label{eq::definitonShapeDerivative}
	d_S \bar G(\kappa) := \frac{d\bar G(\kappa)}{d\kappa},
\end{align}
where $\bar G(\kappa) := \hat G(\Omega_2(\kappa))$.
In the literature one finds formulas for the shape derivative for arbitrary space dimension $d$. In its volume form the shape derivative for some shape $\Omega$ can be written as \cite{LaurainSturm2016}
\begin{align*}
	d_S \hat G(\Omega)(V) = \int_D \mathcal S_1^{\Omega} : \partial V + \mathcal S_0^\Omega \cdot V \ddx,
\end{align*}
where the domain is perturbed by the vector field $V\in C^1_0(D,\R^d)$ with compact support in $D$.
For the considered Poisson problem generalized to arbitrary space dimensions, $\mathcal S_1^\Omega$ and $\mathcal S_0^\Omega$ are given by \cite{LaurainSturm2016}
\begin{align*}
	\mathcal S_1^\Omega =& ((u-\targetU)^2 + \lambda_\Omega \nabla u \cdot \nabla p - f_\Omega p) I - \lambda_\Omega \nabla u \otimes \nabla p - \lambda_\Omega \nabla p \otimes \nabla u, \\
	\mathcal S_0^\Omega =& -2  (u-\targetU) \nabla \targetU - p \nabla f,
\end{align*}
where $p \in H^1_0(D)$ is the solution to the adjoint equation
\begin{align}\label{eq::continousAdjoint}
	\int_D \lambda_\Omega \nabla v \cdot \nabla p \ddx = - 2 \int_D (u - \hat u)v \ddx \quad \mbox{for all } v \in H^1_0(D).
\end{align}
Under sufficient smoothness assumptions, it can be transformed into the Hadamard or boundary form
\begin{align*} 
d_S \hat G(\Omega)(V) = \int_{\overline{\Omega}_1 \cap \overline{\Omega}_2} L \, (V\cdot n) \; \mbox dS_x,
\end{align*}
with $L =( \left. \mathcal S_1^{\Omega} \right \rvert_{\Omega_1} - \left. \mathcal S_1^{\Omega} \right\rvert_{\Omega_2}) n \cdot n$ given by
\begin{align*}
L :=& (\lambda_1- \lambda_2)  (\nabla u \cdot \tau)(\nabla p \cdot\tau) -  \left(\frac{1}{\lambda_1} - \frac{1}{\lambda_2}\right)( \lambda_\Omega \nabla u \cdot n)(  \lambda_\Omega \nabla p \cdot n) . \label{eq_Llambda_inout}
\end{align*}
Here, $n$ denotes the unit normal vector pointing out of $\Omega_1$. It is important to note that the flux $\lambda_\Omega \nabla u \cdot n$ and $\lambda_\Omega \nabla p \cdot n$ are continuous at the interface $\overline \Omega_1 \cap \overline \Omega_2$.
For the 1D setting considered in the present paper we have
\begin{equation}\label{eq::continousShapeDerivativeVolume}
        d_S \bar G(\kappa)(V) = \int_0^{\ell} \left((u-\targetU)^2- fp-\lambda_\Omega \frac{du}{dx} \frac{dp}{dx}\right) \frac{dV}{dx} - \left(2  (u-\targetU)\frac{\targetU}{dx} + \frac{df}{dx}p\right)V  \ddx,
\end{equation}
for a smooth vector field $V$. Choosing $V$ such that $V(\kappa)=1$,
the boundary form of the shape derivative reads
\begin{equation*}
	d_S \bar G(\kappa) = \left(\frac{1}{\lambda_2} - \frac{1}{\lambda_1}\right)\left( \lambda_\Omega \frac{du}{dx} \right)(\kappa) \left(  \lambda_\Omega \frac{dp}{dx} \right)(\kappa).
\end{equation*}
For $f(x)= x$ we have the analytic solution of the state

\begin{equation} \label{eq::analyticU}
u(x)=\left\{\begin{array}{cl} -\frac{x\left(\lambda _{1}\left(-{\ell}^3+\ell x^2+\kappa ^3-\kappa x^2\right)+\lambda _{2}\left(\kappa x^2-\kappa ^3\right)\right)}{6\lambda _{1}\left(\ell\lambda _{1}-\kappa \lambda _{1}+\kappa \lambda _{2}\right)} & \text{\ if\ \ }x\leq \kappa, \\ -\frac{\lambda _{1}\left(\ell-x\right)\left({\ell}^2\kappa -{\ell}^2x+\ell\kappa x-\ell x^2-\kappa ^3+\kappa x^2\right)-\lambda _{2}\left(\ell-x\right)\left({\ell}^2\kappa +\ell\kappa x-\kappa ^3+\kappa x^2\right)}{6\lambda _{2}\left(\ell\lambda _{1}-\kappa \lambda _{1}+\kappa \lambda _{2}\right)} & \text{\ if\ \ }\kappa <x, \end{array}\right.\end{equation}

and for the flux we obtain
\begin{equation}
	q(x) = -\lambda_\Omega \frac{du}{dx} =  \frac{x^2}{2}-\frac{{\ell}^3\lambda _{1}-\kappa ^3\lambda _{1}+\kappa ^3\lambda _{2}}{6\left(\ell\lambda _{1}-\kappa \lambda _{1}+\kappa \lambda _{2}\right)}.
\end{equation}
With \eqref{eq::analyticU} it is possible to obtain explicitly $\bar G(\kappa)$ as a function of $\kappa$, see Appendix \ref{appandix::analyticG}. This allows to compute the shape derivative $d_S \bar G(\kappa)$ symbolically by means of \eqref{eq::definitonShapeDerivative}.
\subsection{Continuous topological derivative}
In this section, we introduce the topological derivative of 
\Cref{prob::tracking}.
Similar as in \Cref{sec::ContinuousShapeDerivative} we focus on a special situation for which we can explicitly compute the topological derivative. In particular, we consider for some $\kappa\in(0,\ell)$ the two-interface situation with interfaces at $\kappa-\pertubationTopology$ and $\kappa+\pertubationTopology$, \ie  $\Omega_2(\kappa,\pertubationTopology) = (\kappa-\pertubationTopology,\kappa+\pertubationTopology)$ and $\Omega_1(\kappa,\pertubationTopology) = (0,\kappa-\pertubationTopology) \cup (\kappa+\pertubationTopology,\ell)$. Such a situation is illustrated in \Cref{fig:flowTwoInterface}. Thus, the continuous topological derivative we investigate in the present paper is defined  by
\begin{align}
	d_T \bar G(\kappa) = \lim_{\pertubationTopology\searrow0}\frac{ \redObjectiveFunction(\Omega_2(\kappa,\pertubationTopology)) - \redObjectiveFunction(\emptyset) }{2\pertubationTopology}.
\end{align}
In literature the general formula for arbitrary space dimension $d$ (with background material $\alpha$ and inclusion material $\beta$, cf. \cite{amstutz2022introduction})
\begin{align*}
	G_\eps(u_\eps) - G_0(u) = & \eps^d \alpha \nabla u(z)^T \mathcal P \nabla p(z) + o(\eps^d)
	\end{align*}
can be found. Here, $\mathcal P$ is the polarization tensor, and $p \in H^1_0(D)$ is the solution to the adjoint problem \eqref{eq::continousAdjoint}.
For the present case, \ie for $d=1$, it holds $\mathcal P = 2(1-\frac{\alpha}{\beta})$ and we have
\begin{align*}
	G_\eps(u_\eps) - G_0(u) =& \overset{|\omega_\eps|}{\overbrace{2\eps}} \alpha \frac{(\beta-\alpha)}{\beta} \frac{du}{dx}(z) \frac{dp}{dx}(z) + o(\eps).
\end{align*}
Therefore, the topological derivative is given by
\begin{align} \label{eq_ana_TD}
	d_T\bar G(\kappa) =& \lambda_1 \frac{\lambda_2-\lambda_1}{\lambda_2} \frac{du}{dx}(\kappa) \frac{dp}{dx}(\kappa).
\end{align}
For $f(x) = x$ and $\hat u(x) = x(\ell-x)$ we have the analytic solution
\begin{equation*}
	d_T\bar G(\kappa) = \bar G_{T1} + \bar G_{T2} + \frac{\bar G_{T3}{\left(\ell-\kappa \right)}^2}{1080\ell{\lambda _{1}}^2\lambda _{2}}
\end{equation*}
with
\begin{align*}
	\bar G_{T1} =& \frac{\kappa ^3\left({\ell}^2-3\kappa ^2\right)\left(\lambda _{1}-\lambda _{2}\right)\left(-10{\ell}^2+60\lambda _{1}\ell+6\kappa ^2-45\lambda _{1}\kappa \right)}{540\ell{\lambda _{1}}^2\lambda _{2}},\\
	\bar G_{T2} =& \frac{\kappa ^2{\left(\ell-\kappa \right)}^2{\left(\ell+\kappa -6\lambda _{1}\right)}^2}{36{\lambda _{1}}^2}, \\
	\bar G_{T3} =& \lambda _{1}\left(7{\ell}^5+14{\ell}^4\kappa -30{\ell}^3\kappa ^2-54{\ell}^2\kappa ^3+27\ell\kappa ^4+36\kappa ^5\right) \\ & -\lambda _{2}\left(7{\ell}^5+14{\ell}^4\kappa +6{\ell}^2\kappa ^3+57\ell\kappa ^4+36\kappa ^5\right) \\ & -{\lambda _{1}}^2\left(30{\ell}^4+60{\ell}^3\kappa -180{\ell}^2\kappa ^2-180\ell\kappa ^3+270\kappa ^4\right) \\ &+\lambda _{1}\lambda _{2}\left(30{\ell}^4+60{\ell}^3\kappa +180{\ell}^2\kappa ^2+180\ell\kappa ^3+270\kappa ^4\right) \\ &-1080\ell\kappa ^2{\lambda _{1}}^2\lambda _{2}.
\end{align*}
\section{Discretization}\label{sec::discretization}
\begin{figure}
	\begin{subfigure}{0.33\textwidth}
		\centering
		\includegraphics[width=0.9\linewidth]{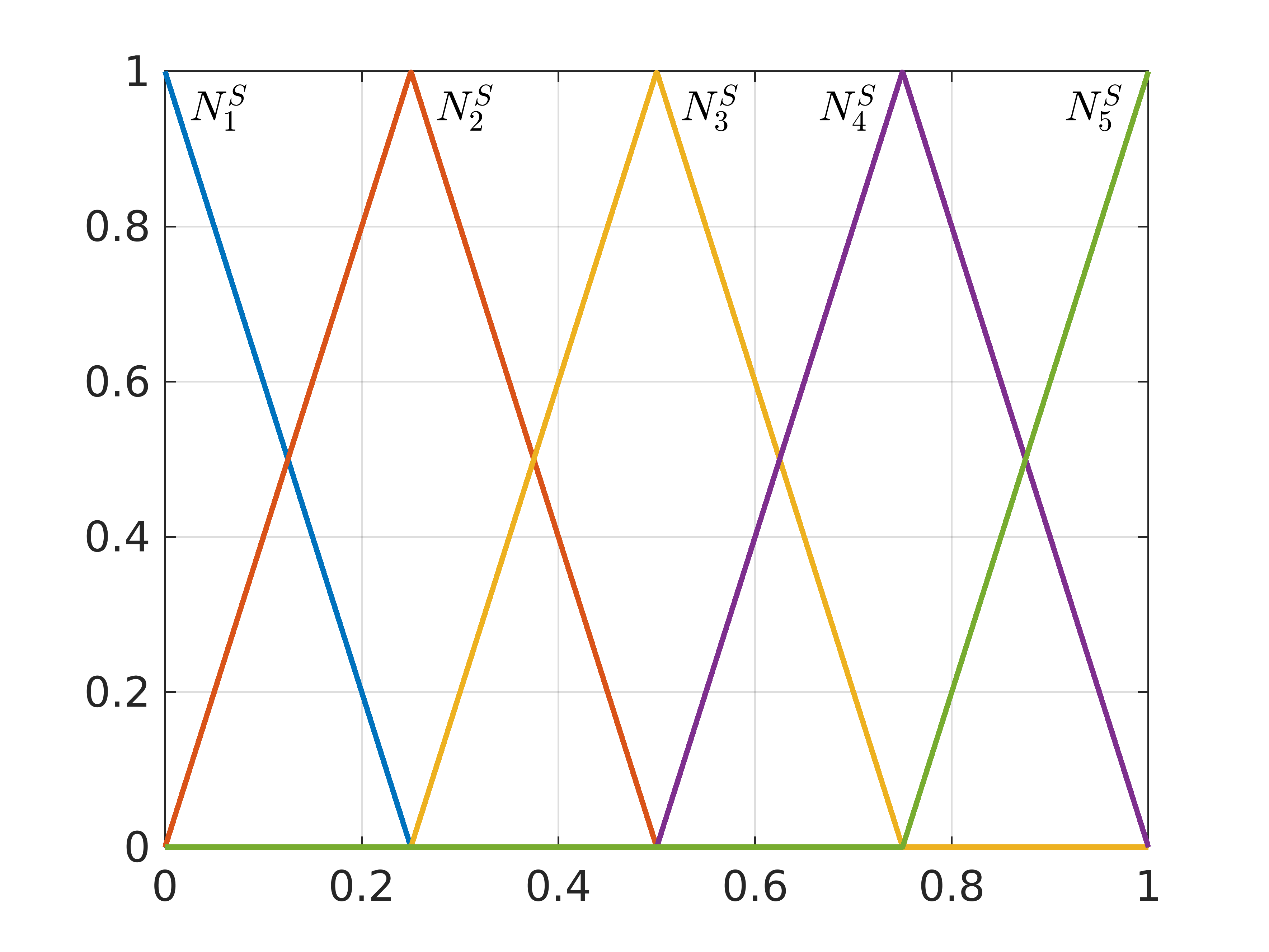}
		\caption{$p=1$}
	\end{subfigure}\hfil
	\begin{subfigure}{0.33\textwidth}
		\centering
		\includegraphics[width=0.9\linewidth]{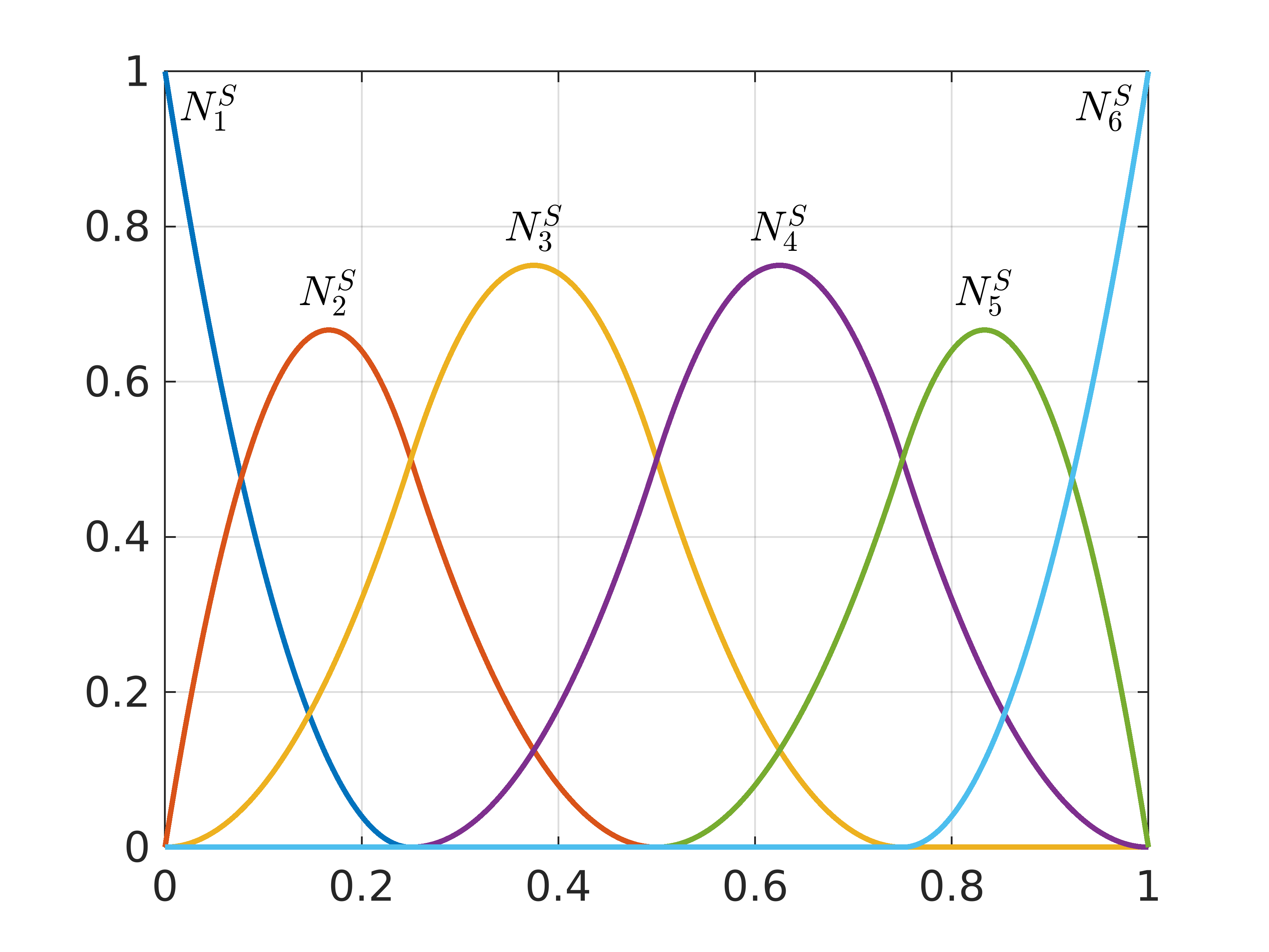}
		\caption{$p=2$}
	\end{subfigure}\hfil
	\begin{subfigure}{0.33\textwidth}
		\centering
		\includegraphics[width=0.9\linewidth]{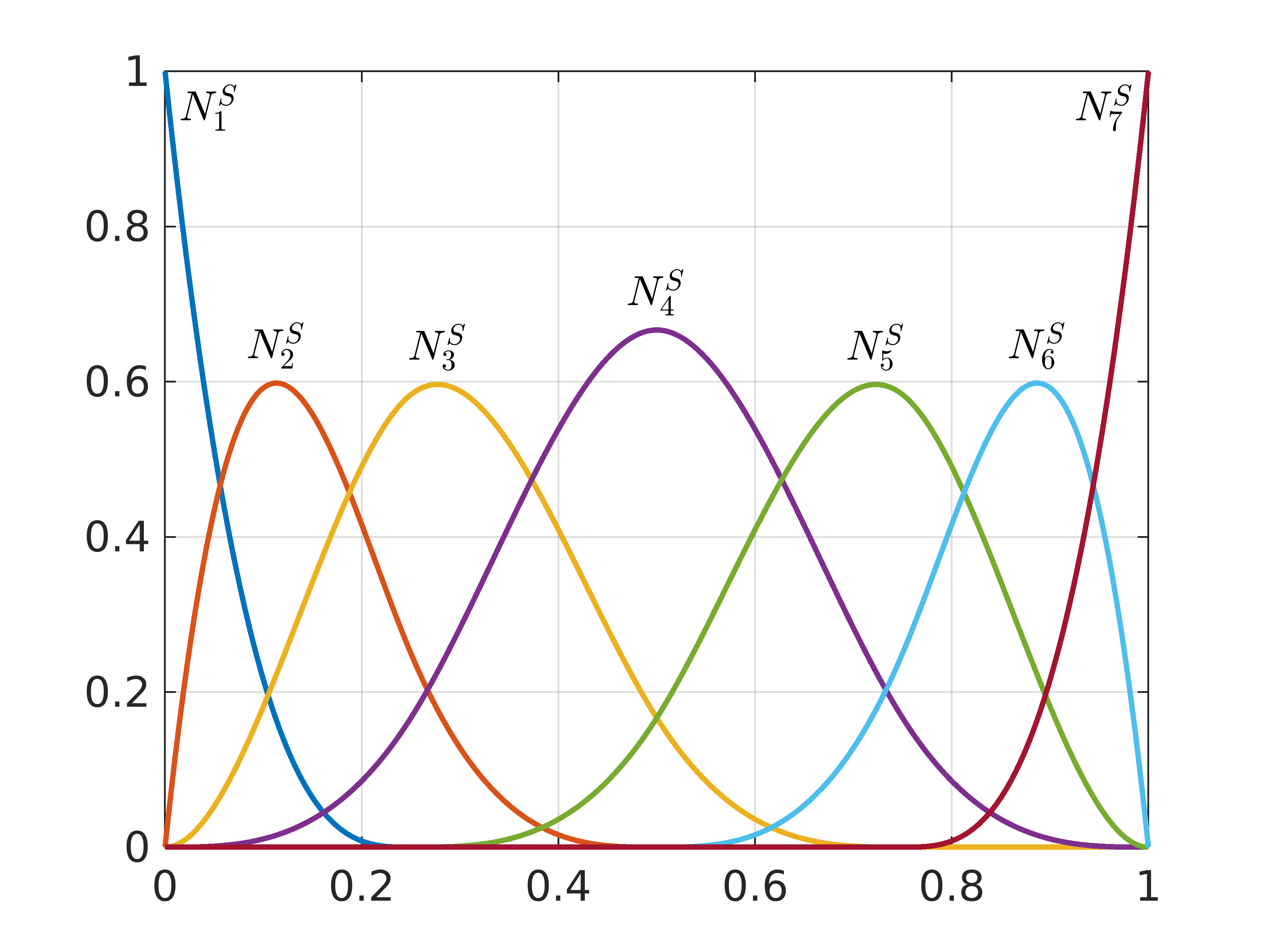}
		\caption{$p=3$}
	\end{subfigure}
	\caption{Discretization without enrichment functions}
	\label{fig:StandardDiscretizations}
\end{figure}

\Cref{prob::tracking} is discretized by B-splines with polynomial  degree $p=1,2,3$ and smoothness $C^{p-1}$. For fixed p, the $n_S$ basis functions $\varphi_{i,p}(\xi)$, $i=1,\dots,n_S$ are defined over the parametric domain which is defined by a knot vector. The knot vector is a non-decreasing
set of real numbers called knots $\pmb\zeta = \{\zeta_1,\zeta_2,\dots,\zeta_{n_S+p+1}\}$. Here, we use an open knot vector where the first and last knots are repeated p + 1 times. For further details on B-splines, see e.g.\cite{1995piegla}. In the \standardMethod the approximation reads
\begin{equation}
	u^S(x) = \sum_{i=1}^{n_S} N^S_i(x) u^S_i,
\end{equation}
with the basis functions $N^S_i$ of degree $p$ implicitly defined by $N^S_i(x(\xi)) = \varphi_{i,p}(\xi)$ (see \Cref{fig:StandardDiscretizations}).
\begin{figure}
	\begin{subfigure}{0.33\textwidth}
		\centering
		\includegraphics[width=0.9\linewidth]{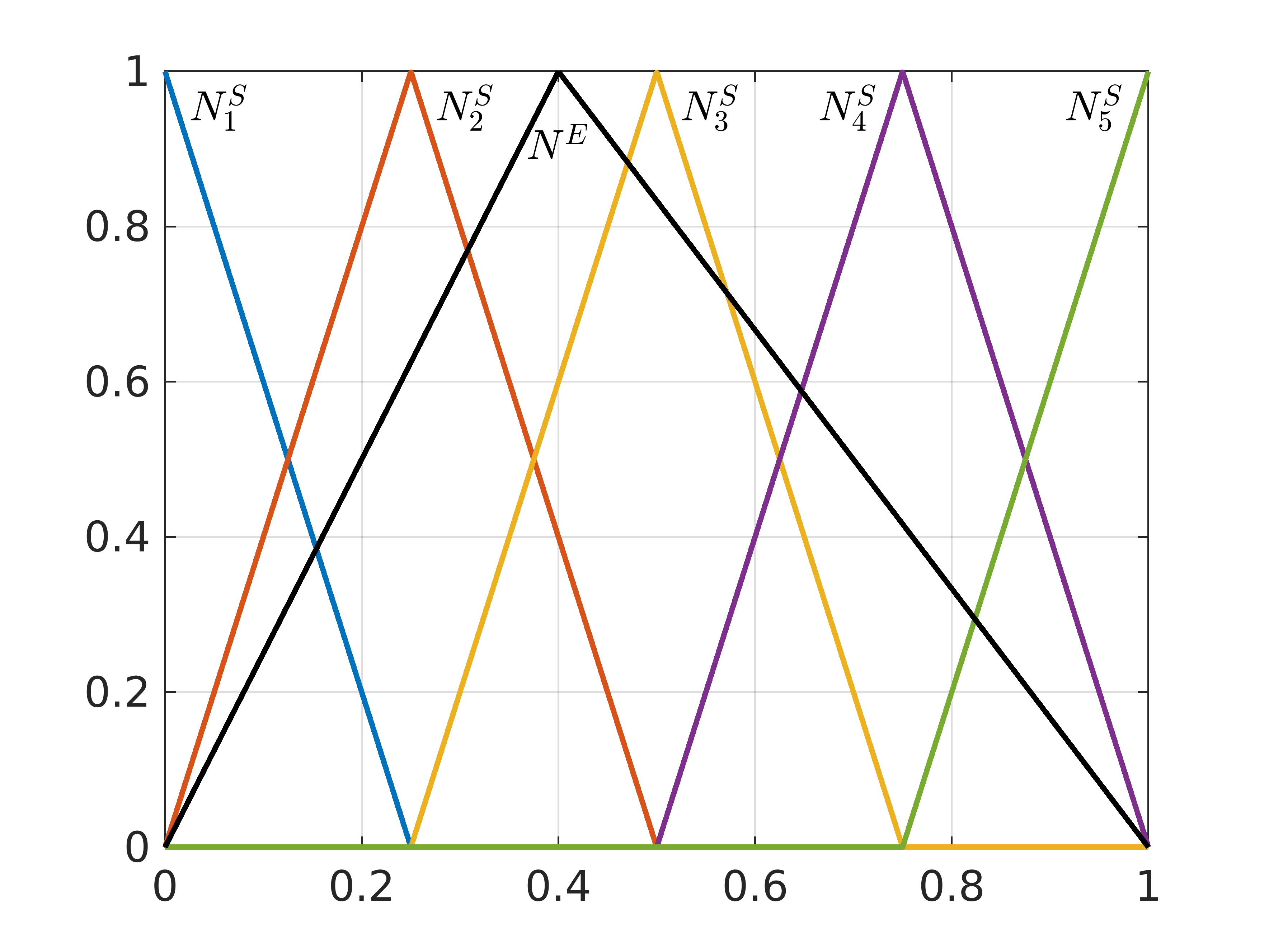}
		\caption{$p=1$}
	\end{subfigure}\hfil
	\begin{subfigure}{0.33\textwidth}
		\centering
		\includegraphics[width=0.9\linewidth]{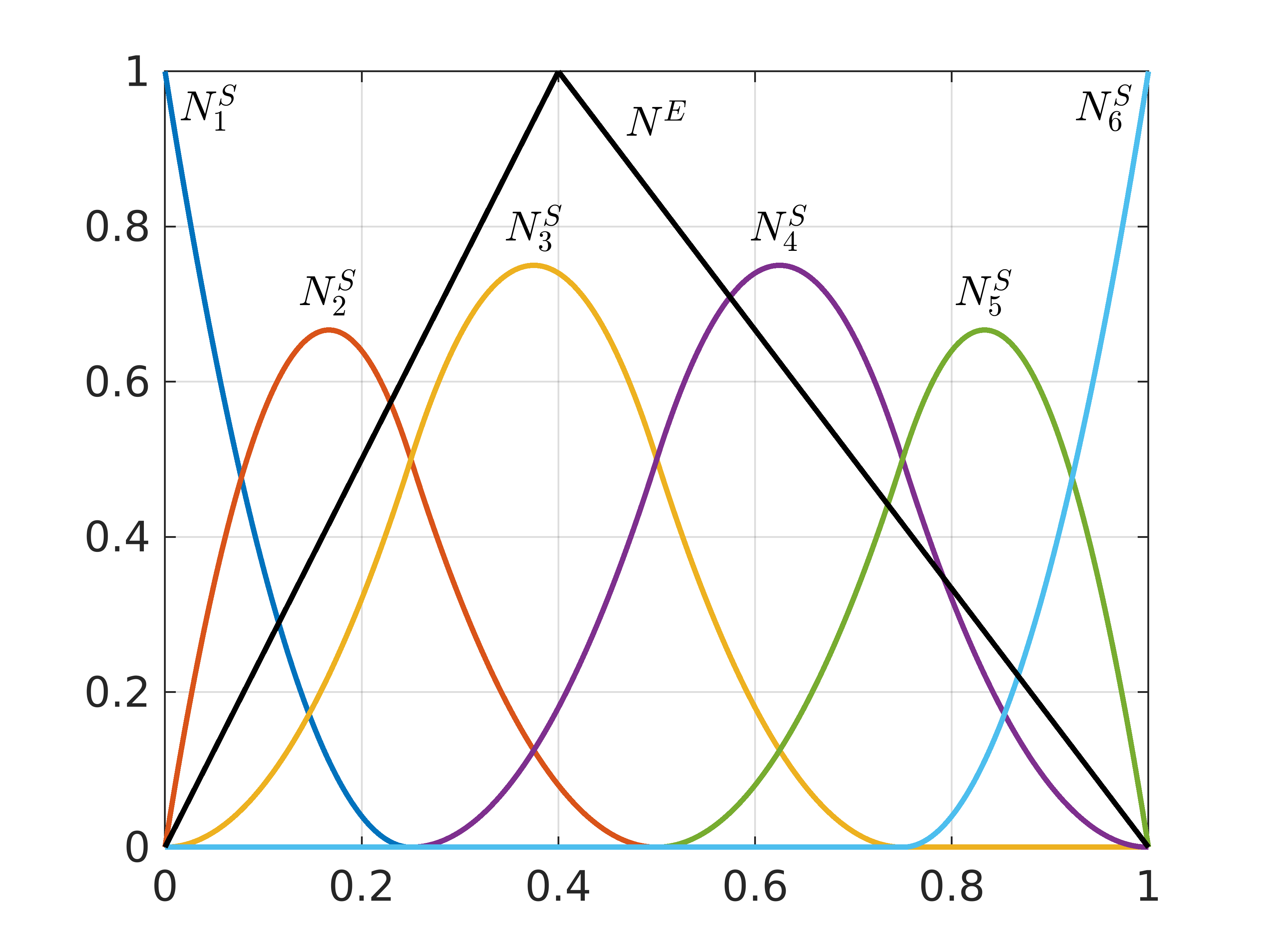}
		\caption{$p=2$}
	\end{subfigure}\hfil
	\begin{subfigure}{0.33\textwidth}
		\centering
		\includegraphics[width=0.9\linewidth]{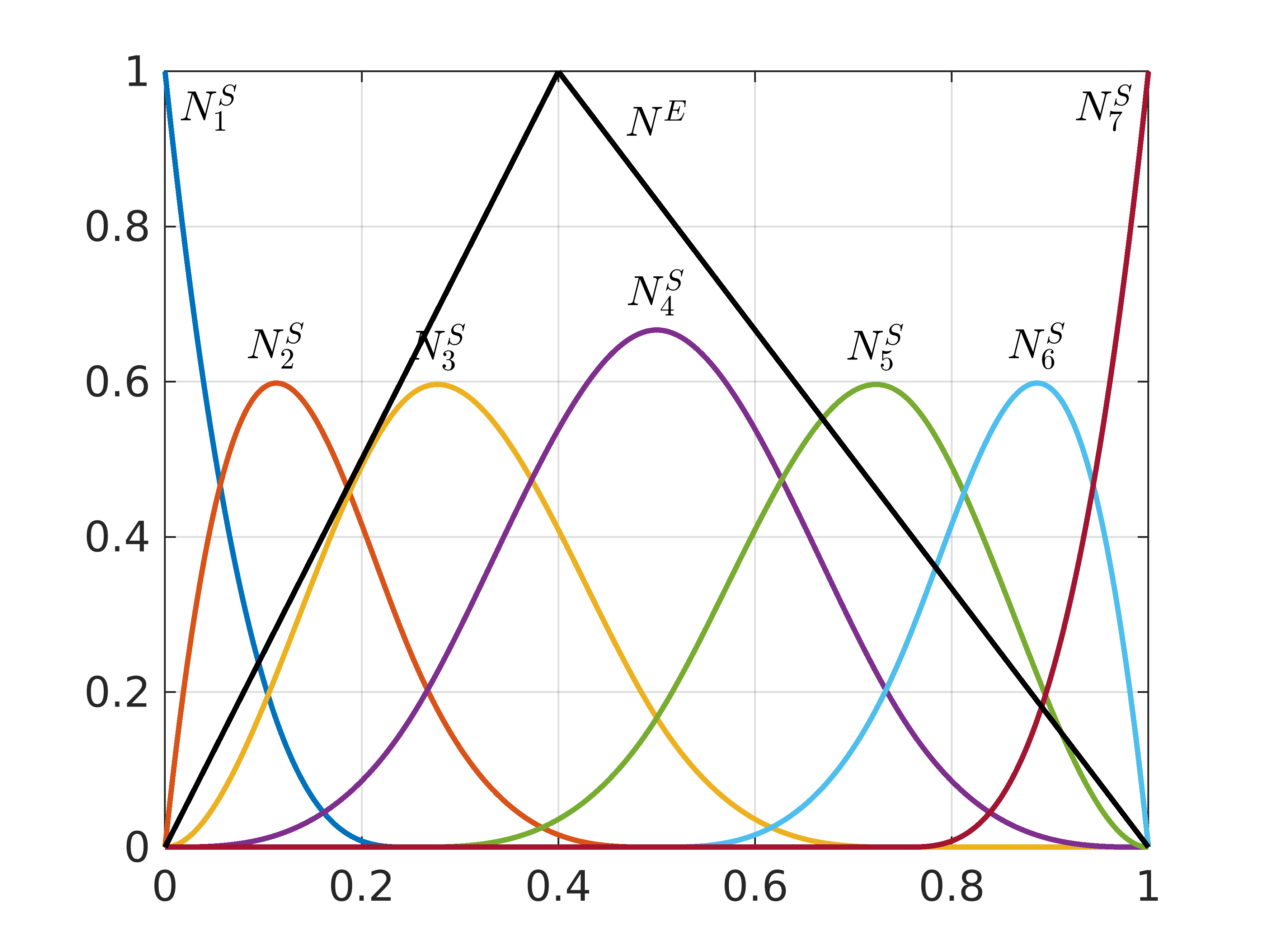}
		\caption{$p=3$}
	\end{subfigure}
	\caption{Enriched discretization with one interface at $x=0.4$. The enrichment function is denoted by $N^E$.}
	\label{fig:EnrichedDiscretizations}
\end{figure}
\corrections{For the \enrichedMethod we consider the set of interface points $\overline \Omega_1\cap \overline \Omega_2=\{\kappa_1\}$ in the one-interface case and set $n_E:=1$, or $\overline \Omega_1\cap \overline \Omega_2=\{\kappa_1, \kappa_2\}$ in the two-interface case where $n_E:=2$. Here,}
we use a global linear enrichment strategy such that the approximation reads
\begin{equation}\label{eq:standard}
	u^E(x) = \sum_{i=1}^{n_S} N^S_i(x)  u^S_i + \sum_{j=1}^{n_E}\enrichmentBasis_j(x)  u^E_j,
\end{equation}
where the enrichment function $\enrichmentBasis_j$, $j=1,\dots,n_E$ is associated with the interface at $\kappa_j \in \overline \Omega_1\cap \overline \Omega_2$,
\begin{equation}\label{eq:enrichment}
	\enrichmentBasis_j(x) = \begin{cases}
		\frac{x}{\kappa_j}  &\text{ for } x \le \kappa_j, \\
		\frac{\ell-x}{(\ell-\kappa_j)} &\text{ for } x > \kappa_j.
	\end{cases}
\end{equation}
We refer to \Cref{fig:EnrichedDiscretizations} for an illustration of the basis functions for the \enrichedMethod. \corrections{Note that in \eqref{eq:standard} and \eqref{eq:enrichment} the same basis functions $N^S_i$ are used. In \eqref{eq:enrichment}, one global enrichment function $\enrichmentBasis_j$ per interface point $\kappa_j$ is added.}

In order to treat the two methods in a unified way we denote the $n = n_S + n_E$ basis functions by $\psi_i$, which span the approximation space $V_h(D)$. More precisely, we have $\psi_i = N_i^S$ for $i=1,\dots, n_S$ and $\psi_i = N_{i-n_S}^E$ for $i = n_S+1, \dots, n$. For the standard method we have $n_E=0$.
Next, we give the discretized state problem.
\begin{mdframed}
\begin{problem}\label{prob::discretized}
	For fixed $\ell,\lambda_1, \lambda_2>0$, $f:D\rightarrow\R$ the discretized state problem reads:
    \begin{align*}
	&\mbox{Find }u_h\in V_h(D) \mbox{ such that } \nonumber \\
	&\int_D \lambda_\Omega \frac{d u_h}{dx}(x)  \frac{d v_h}{dx}(x) \ddx = \int_D f(x)\, v_h(x) \ddx \quad  \mbox{ for all }v_h \in V_h(D).
	\end{align*}
\end{problem}
\end{mdframed}
The system of linear equations corresponding to \Cref{prob::discretized} reads
\begin{equation*}
    \mathbf K \mathbf u = \mathbf f,
\end{equation*}
with the solution vector $\mathbf u \in \R^{n}$, the stiffness matrix $\mathbf K \in \R^{n\times n}$, and the right-hand-side vector $\mathbf f \in \R^{n}$ given by
\begin{equation*}
	\begin{aligned}
        \mathbf K[i,j] &= \int_{D} \lambda_{\Omega} \frac{d\psi_j}{dx}(x)   \frac{d\psi_i}{dx}(x) \ddx,  \\
		\mathbf f[i] &= \int_{D}  f(x)\psi_i(x) \,   \ddx.
		\end{aligned}
\end{equation*}
For the statement of the discretized optimization problem we define the mass matrix $\mathbf M \in \R^{n\times n}$ and the vector $\mathbf m \in \R^{n}$, given by
\begin{equation*}
	\begin{aligned}
		\mathbf M[i,j] &= \int_{D} \psi_j(x)  \psi_i(x) \ddx \\
        \mathbf m[i]  &= \int_{D}  \hat u(x) \psi_i(x) \ddx, \\
		\end{aligned}
\end{equation*}
Next, the discretized optimization problem is given.
\begin{mdframed}
\begin{problem}\label{prob::trackingDiscrete}
	For fixed parameters $\ell,\lambda_1, \lambda_2>0$, the source term $f:D\rightarrow\R$ and the desired state $\targetU:D\rightarrow\R$ the discretized optimization problem reads:
		\begin{alignat*}{2}
			\min_{\Omega_2\in \mathcal A} \tilde\objectiveFunction(\mathbf u) &= \mathbf u^\top\MM\mathbf u - 2 \mathbf u^\top \mathbf m \\
			\text{ subject to}& \nonumber \\
			\mathbf u \in \R^n \; &\mbox{is the coefficient vector of the solution of} \; \Cref{prob::discretized}.\nonumber
		\end{alignat*}
\end{problem}
\end{mdframed}

\section{Sensitivity analysis of the discretized problem}
\label{sec::numericalSensitivities}
In this section, we perform shape and topological sensitivity analysis of \Cref{prob::trackingDiscrete}.
To this end, we use the Lagrangian approach of \cite{gangl2023unified} in \Cref{sec::discreteShapeStandard,sec::discreteShapeEnriched,sec::discreteTopStandard}. However, for the topological sensitivity analysis of the \enrichedMethod in \Cref{sec::discreteTopEnriched} we develop a new approach based on the direct method.
Let $\mathbf K_\eps$, $\MM_\eps$, $\ff_\eps$ and $\mm_\eps$ be the stiffness matrix, the mass matrix, the right-hand-side vector, and the vector involving $\targetU$, respectively, after a shape or topological perturbation. Note that in the case of the standard method, geometric perturbations only affect the stiffness matrix via the diffusion coefficient. On the other hand, in the case of the enriched method also the enrichment functions $\psi_{n_S+1}, \dots, \psi_{n_S+n_E}$ depend on the position of the interface. Thus $\MM_\eps \neq \MM$ and similar for $\ff_\eps$ and $\mm_\eps$ in this case. Let $\mathbf u_\eps$ be the solution of
\begin{equation*}
    \mathbf K_\eps \mathbf u_\eps = \mathbf f_\eps,
\end{equation*}
and $\pp$ the adjoint solution to \Cref{prob::trackingDiscrete} defined by
\begin{equation*}
    \mathbf K^T \pp =-2(\MM\mathbf u-\mm).
\end{equation*}

By adapting the results of \cite{gangl2023unified} for the numerical \topshapeDerivative derivative, we obtain for the \standardMethod in the cases of shape and topological perturbations, as well as for the \enrichedMethod in the case of a shape perturbation, the formula
\begin{align}\label{eq::derivativeTracking}
	d \bar G(\kappa) &=
	\frac{\mathbf p^\top(d \Amat \,\uuzero - d \mathbf f)
		+ \uuzero ^\top d\MM  \uuzero - 2 \uuzero ^\top d\mathbf{m}	}{ d a},
\end{align}
where for the 1D problem setting considered in the present paper%
\begin{align} \label{eq_derivativeLimits}
	\begin{aligned}
		d \Amat &= \lim_{\varepsilon\searrow 0} \frac{\Amat_\eps- \Amat}{\varepsilon}, &&&
		d \MM& = \lim_{\varepsilon\searrow 0} \frac{\MM_\eps- \MM}{\varepsilon}, \\
		d \ff &= \lim_{\varepsilon\searrow 0} \frac{\ff_\eps-\ff}{\varepsilon}, &&&
		d\mathbf m &= \lim_{\varepsilon\searrow 0} \frac{\mathbf m_\eps- \mathbf m}{\varepsilon},&&&   \\
    & &&&d  a &= \begin{cases}
		    1 \quad \mbox{for the shape derivative,}\\
            2 \quad \mbox{for the topological derivative.}
		\end{cases}
	\end{aligned}
\end{align}

\begin{remark}
    Note that in \cite{gangl2023unified} it is assumed that the desired state $\targetU$ is discretized and represented by the vector $\hat{\uu}$. In contrast to this setting of \cite{gangl2023unified}, we treat $\targetU$ here as given continuous function and have incorporated it in \Cref{prob::trackingDiscrete} by means of the vector $\mathbf m$.
\end{remark}

\subsection{Numerical shape derivative}\label{sec::discreteShape}
Next, we compute the numerical \topshapeDerivative derivative given in \eqref{eq::derivativeTracking} for the case of a shape perturbation for both the \standardMethod and the \enrichedMethod. By a shape perturbation, we mean that the interface is perturbed a distance $\eps$ from its original position $\kappa$ into a direction $\delta \in \{-1,1\}$, i.e., to a perturbed position $\kappa + \delta \eps$.
For notational simplicity, we drop the direction $\delta$ and just consider the right and left sided limit with respect to $\eps$.

\subsubsection{Shape derivative for the problem discretized by the \standardMethod}\label{sec::discreteShapeStandard}
For the \standardMethod we have
	\begin{align*}
		\Amat_\eps = \Amat^S_\eps \in \R^{n_S\times n_S}
	\end{align*}
with
\begin{align}\label{eq::KSeps}
	\mathbf K^S_\eps[i,j] &=
	\int_0^{\kappa+\eps} \lambda_1 \frac{dN^S_i(x)}{dx} \frac{dN^S_j(x)}{dx} \ddx
    + \int_{\kappa+\eps}^\ell \lambda_2 \frac{dN^S_i(x)}{dx} \frac{dN^S_j(x)}{dx} \ddx.
\end{align}
For smooth basis functions (i.e. $p>1$) or if the interface does not coincide with a node (i.e. $\kappa \neq x_k, \; 1,...,n_S$) we obtain by the mean value theorem for integrals
\begin{align*}
		d\mathbf K[i,j] = \lim_{\eps\searrow0}\frac{\mathbf K^S_\eps - \mathbf K^S_0}{\eps}[i,j] &= \lim_{\eps\searrow0}\frac{\lambda_1-\lambda_2}{\eps}\int_\kappa^{\kappa+\eps}  \frac{dN^S_i(x)}{dx} \frac{dN^S_j(x)}{dx} \ddx \\
		&=(\lambda_1-\lambda_2)\frac{dN^S_i}{dx}(\kappa)  \frac{dN^S_j}{dx}(\kappa)
\end{align*}
for the right-sided limit. In this case, the left sided limit is easily checked to coincide up to the sign.
All other terms in \eqref{eq::derivativeTracking} are zero and we obtain
\begin{align}\label{eq::derivativeTrackingShapeStandard}
	d_S \bar G(\kappa) &= \mathbf p^\top d \Amat \,\uuzero.
\end{align}
\begin{remark}\label{remark::shape}
    For $p=1$ the shape derivative is not continuous if the interface coincides with a node. In this case, one can compute the left- and the right-sided limits
\begin{align*}
    d\mathbf K^+[i,j] = (\lambda_1-\lambda_2)\frac{dN^S_i}{dx}(\kappa^+)  \frac{dN^S_j}{dx}(\kappa^+), \\
    d\mathbf K^-[i,j] = (\lambda_1-\lambda_2)\frac{dN^S_i}{dx}(\kappa^-)  \frac{dN^S_j}{dx}(\kappa^-).
\end{align*}
\end{remark}
\subsubsection{Shape derivative for the problem discretized by the \enrichedMethod}\label{sec::discreteShapeEnriched}
For the \enrichedMethod the stiffness matrix has the block structure
	\begin{align}
		\Amat_\eps =
		\begin{bmatrix}
			\KmatSeps & \KmatSEeps \\
			(\KmatSEeps)^\top & \KEeps
		\end{bmatrix}\in \R^{(n_S+n_E)\times (n_S+n_E)},
	\end{align}
with $\KmatSeps$ given by \eqref{eq::KSeps}. The entries of the additional blocks are given by
	\begin{align*}
		\mathbf K^{SE}_\eps[i] &=
		\int_0^{\kappa+\eps} \lambda_1 \frac{dN^S_i(x)}{dx} \frac{1}{\kappa+\eps} \ddx + \int_{\kappa+\eps}^\ell \lambda_2 \frac{dN^S_i(x)}{dx}\left(-\frac{1}{\ell-\kappa-\eps}\right)  \ddx, \\
		  \KEeps &=
		\int_0^{\kappa+\eps} \lambda_1 \left(\frac{1}{\kappa+\eps}\right)^2 \ddx +
		\int_{\kappa+\eps}^\ell \lambda_2\left(-\frac{1}{\ell-\kappa-\eps}\right)^2 \ddx.
	\end{align*}
    Thus, \Cref{remark::shape} applies also for this section and in particular also for the entries of the matrix $\KmatSEeps$.
	Due to the assumption of only one interface we have $n_E=1$ and thus $\KEeps \in \R$ and $\KmatSEeps \in \R^{n_S \times 1}$. We proceed by computing the derivatives of the stiffness matrix for $d \Amat$. We have
	\begin{equation*}
	\begin{aligned}
		(\KmatSEeps - \KmatSE)[i] &= \int_0^{\kappa} \lambda_1 \frac{dN^S_i(x)}{dx} \left(\frac{1}{\kappa+\eps}-\frac{1}{\kappa}\right) \ddx  \\
		&+ \int_\kappa^{\kappa+\eps} \frac{dN^S_i(x)}{dx}\left(\frac{\lambda_1 }{\kappa+\eps}+\frac{\lambda_2 }{\ell-\kappa}\right) \ddx \\
		& - \int_{\kappa+\eps}^\ell \lambda_2\frac{dN^S_i(x)}{dx}\left(\frac{1}{\ell-\kappa-\eps}-\frac{1}{\ell-\kappa}\right) \ddx,
	\end{aligned}
	\end{equation*}
	and considering
		\begin{align*}
			\frac{1}{\kappa+\eps}-\frac{1}{\kappa} &= \frac{\kappa-\kappa-\eps}{(\kappa+\eps)\kappa} = \frac{-\eps}{(\kappa+\eps)\kappa}, \\
			\frac{1}{\ell-\kappa-\eps}-\frac{1}{\ell-\kappa} &= \frac{(\ell-\kappa)-(\ell-\kappa-\eps)}{(\ell-\kappa-\eps)(\ell-\kappa)} = \frac{\eps}{(\ell-\kappa-\eps)(\ell-\kappa)},
		\end{align*}
	yields
	\begin{align}
	\begin{aligned}
		\lim_{\eps\searrow0}\frac{\KmatSEeps - \KmatSE}{\eps}[i]&= -\frac{\lambda_1}{\kappa^2}\int_0^{\kappa}\frac{dN^S_i(x)}{dx}\ddx-\frac{\lambda_2 }{(\ell-\kappa)^2}\int_\kappa^{\ell}\frac{dN^S_i(x)}{dx}\ddx + \frac{dN^S_i}{dx}(\kappa^+)\left(\frac{\lambda_1}{\kappa}+\frac{\lambda_2}{\ell-\kappa}\right) \\
        &=\lambda_1\int_0^{\kappa}\frac{dN^{E}(x)}{dxd\kappa}\frac{dN^S_i(x)}{dx}\ddx + \lambda_2  \int_\kappa^{\ell} \frac{dN^{E}(x)}{dxd\kappa}\frac{dN^S_i(x)}{dx}\ddx \\
        &\quad+ \lambda_1 \frac{dN^S_i}{dx}(\kappa^+) \frac{dN^{E}}{dx}(\kappa^-) + \lambda_2 \frac{dN^S_i}{dx}(\kappa^+) \frac{dN^{E}}{dx}(\kappa^+).
	\end{aligned}
	\end{align}
	Furthermore,
	\begin{align*}
		\KEeps - \KE=& \int_0^{\kappa} \lambda_1  \left(\frac{1}{(\kappa+\eps)^2}-\frac{1}{\kappa^2}\right) \ddx  \\
		&+ \int_\kappa^{\kappa+\eps} \left(\frac{\lambda_1 }{(\kappa+\eps)^2}-\frac{\lambda_2 }{(\ell-\kappa)^2}\right) \ddx \\
		& - \int_{\kappa+\eps}^\ell \lambda_2\left(\frac{1}{(\ell-\kappa-\eps)^2}-\frac{1}{(\ell-\kappa)^2}\right) \ddx,
	\end{align*}
	and with
		\begin{align*}
			\frac{1}{(\kappa+\eps)^2}-\frac{1}{\kappa^2} &= \frac{\kappa-\kappa-\eps}{(\kappa+\eps)\kappa} = \frac{-2\kappa\eps-\eps^2}{(\kappa+\eps)^2\kappa^2}, \\
			\frac{1}{(\ell-\kappa-\eps)^2}-\frac{1}{(\ell-\kappa)^2} &= \frac{2(\ell-\kappa)\eps-\eps^2}{(\ell-\kappa-\eps)^2(\ell-\kappa)^2},
		\end{align*}
	the limit yields
	\begin{align}
		\lim_{\eps\searrow0}\frac{\KEeps - \KE}{\eps}&= -\frac{\lambda_1}{\kappa^2}+\frac{\lambda_2 }{(\ell-\kappa)^2}.
	\end{align}
Thus, for the \enrichedMethod even for $\lambda_1=\lambda_2$ we have a non-zero shape derivative. Furthermore, although $f(x)$ is independent of $\kappa$ in \Cref{prob::tracking}, we have to consider also the vector $d\mathbf f = [d\mathbf f^S,df^E]^\top$. Here, $d\mathbf f^S=\mathbf 0$ but for $ df^E$ we obtain
\begin{align}
	\lim_{\eps\searrow0}\frac{ f^E_\eps- f^E_0}{\eps}&= \int_0^{\ell} f(x)\lim_{\eps\searrow0}\frac{\enrichmentBasis_{\kappa+\eps}(x)-\enrichmentBasis_{\kappa}(x)}{\eps} \ddx
	= \int_0^{\ell} f(x)\frac{d\enrichmentBasis_{\kappa}(x)}{d\kappa} \ddx .
\end{align}
For the \enrichedMethod  we also have a block structure for the mass matrix,
\begin{align}
	\MM_\eps =
	\begin{bmatrix}
		\MM^S_\eps & \MM^{SE}_\eps \\
		(\MM^{SE}_\eps)^\top & M^E_\eps
	\end{bmatrix},
\end{align}
with
\begin{subequations}
\begin{align*}
	\MM^S_\eps[i,j] = \MM^S[i,j] &= \int_{0}^{\ell} N^S_i(x) N^S_j(x) \ddx, \\
	\MM^{SE}_\eps[i] &= \int_{0}^{\ell} N^S_i(x) \enrichmentBasis_{\kappa+\eps}(x) \ddx, \\
	M^E_\eps &= \int_{0}^{\ell}  \enrichmentBasis_{\kappa+\eps}(x) \enrichmentBasis_{\kappa+\eps}(x) \ddx = \frac{\ell}{2}.
\end{align*}
\end{subequations}
Furthermore, the vector $\mathbf m_\eps$ has the structure
\begin{align}
	\mathbf m_\eps =
	\begin{bmatrix}
		\mathbf m^S_\eps   \\
		m^E_\eps
	\end{bmatrix},
\end{align}
with
\begin{subequations}
\begin{align*}
	\mathbf m^S_\eps[i] = \mathbf m^S[i] &= \int_{0}^{\ell} N^S_i(x) \hat u(x) \ddx, \\
	m^E_\eps &= \int_{0}^{\ell}  \enrichmentBasis_{\kappa+\eps}(x) \hat u(x) \ddx.
\end{align*}
\end{subequations}
Thus, we obtain
\begin{align}
	d\MM =
	\begin{bmatrix}
		\mathbf 0 & d\MM^{SE}_\eps \\
		(d\MM^{SE}_\eps)^\top & 0
	\end{bmatrix}, \quad
    d\mm = \begin{bmatrix}
		\mathbf 0   \\
		d m^{E}
	\end{bmatrix}
\end{align}
with the limits
\begin{align}
	d\mathbf M^{SE}[i] = \lim_{\eps\searrow 0}\frac{\MmatSEeps-\MmatSE}{\eps}[i]&= \int_{0}^{\ell}  N^S_i(x) \frac{d\enrichmentBasis_{\kappa}(x)}{d\kappa} \ddx, \\
	d m^{E}=\lim_{\eps\searrow 0}\frac{ m^E_\eps- m^E_0}{\eps} &= \int_{0}^{\ell}\frac{d\enrichmentBasis_{\kappa}(x)}{d\kappa} \hat u(x) \ddx.
\end{align}
To sum up, we have computed that the shape derivative for the \enrichedMethod is given by
\begin{equation}\label{eq::derivativeTrackingShapeEnriched}
	d_S \bar G(\kappa) = \mathbf p^\top d \Amat \,\uuzero + df^E p^E + 2 (u^E d\mathbf M^{SE} \mathbf u^S -  dm^E u^E).
\end{equation}

\begin{remark}
    While all these computations were done for the right-sided limit, left sided limit can be seen to coincide by straightforward calculations except for the case mentioned in Remark \ref{remark::shape}.
\end{remark}

\subsection{Numerical topological derivative}\label{sec::discreteTop}
In this section, we consider the case where a homogeneous material is perturbed by a topological perturbation around a position $\kappa \in (0, \ell)$, i.e., in the unperturbed configuration we have $\Omega_1 = (0, \ell)$ and $\Omega_2=\emptyset$, see also \Cref{fig:flowTwoInterface}.

The computation of numerical topological derivatives involves limits of the form
\begin{equation*}
    \lim_{\eps \searrow 0} \frac{u_\eps - u}{\eps},
\end{equation*}
were $u_\eps$ and $u$ are the states for the perturbed and unperturbed problems respectively.
Thus, a necessary condition for these computations is that
\begin{equation} \label{eq_limueps_u}
	\lim_{\eps\searrow 0} u_\eps = u.
\end{equation}

In the case of a discretization by the \enrichedMethod, if $\kappa$ is within an interval $(x_{k-1}, x_k)$, we have in the limit a global linear enrichment function with node at $\kappa$. This function is not in the span of the standard basis and thus \eqref{eq_limueps_u} cannot hold in this case. Moreover, even if $\kappa \in \{x_1, \dots, x_n\}$, for spline degree $p>1$ the linear enrichment function involved in $u_\eps$ cannot be represented by the smooth standard basis functions composing $u$ such that \eqref{eq_limueps_u} cannot hold in this case either.
For these reasons, we restrict ourselves to the particular case $\kappa = x_k$ for some $k \in \{2, \dots, n-1\}$ and spline degree $p=1$ in \Cref{sec::discreteTopEnriched}.

Since we are interested in a comparison between the \standardMethod and the \enrichedMethod, we also restrict ourselves to this setting in the case of the \standardMethod.

\subsubsection{Topological derivative for the \standardMethod}\label{sec::discreteTopStandard}
In this section, we compute the topological derivative for the \standardMethod using the approach from \cite{gangl2023unified}. Assuming a topological change at an interior node $x_k$ the perturbed stiffness matrix reads
\begin{equation}
    \mathbf K_\eps[i,j] = \int_0^{x_k-\eps} \lambda_1 \frac{dN^S_i}{dx} \frac{dN^S_j}{dx} \ddx
    + \int_{x_k-\eps}^{x_k+\eps} \lambda_2 \frac{dN^S_i}{dx} \frac{dN^S_j}{dx} \ddx
    + \int_{x_k+\eps}^{\ell} \lambda_1 \frac{dN^S_i}{dx} \frac{dN^S_j}{dx} \ddx.
\end{equation}
For the difference between perturbed and unperturbed stiffness matrix we obtain
\begin{equation*}
    (\mathbf K_\eps-\mathbf K_0)[i,j] =  (\lambda_2-\lambda_1)\int_{x_k-\eps}^{x_k+\eps}  \frac{dN^S_i}{dx} \frac{dN^S_j}{dx} \ddx,
\end{equation*}
and therefore
\begin{equation}
	\begin{aligned}
		d\mathbf K[i,j] = \lim_{\eps\searrow0}\frac{\mathbf K_\eps - \mathbf K_0}{\eps}[i,j]
		&=\left(\lambda_2-\lambda_1\right)\left(\frac{dN^S_i}{dx}(x_k^-)  \frac{dN^S_j}{dx}(x_k^-) + \frac{dN^S_i}{dx}(x_k^+)  \frac{dN^S_j}{dx}(x_k^+)\right),
	\end{aligned}
\end{equation}
where
 \begin{equation}
    \frac{dN^S_j}{dx}(x_k^-) = \begin{cases}
    \frac{1}{h_k} \quad &\mbox{if } j = k,\\
    -\frac{1}{h_k} \quad &\mbox{if } j = k-1,\\
    0, \quad &\mbox{else}
\end{cases}\quad
\frac{dN^S_j}{dx}(x_k^+) = \begin{cases}
    -\frac{1}{h_{k+1}} \quad &\mbox{if } j = k,\\
    \frac{1}{h_{k+1}} \quad &\mbox{if } j = k+1,\\
    0. \quad &\mbox{else}
\end{cases}
\end{equation}
Here, $h_k = x_k - x_{k-1}$ and $h_{k+1} = x_{k+1} - x_{k}$.
Again, all other terms then $d\mathbf K$ in \eqref{eq::derivativeTracking} are zero for the \standardMethod. Thus, the topological derivative for node $x_k$ is given by
\begin{equation}\label{eq::topologicalDerivativeStandard}
	d_T\redObjectiveFunction(x_k) = \frac{\lambda_2-\lambda_1}{2} \begin{bmatrix}
	    p_{k-1} \\ p_k \\ p_{k+1}
	\end{bmatrix}^\top \begin{bmatrix}
	    h_k^{-2} & -h_k^{-2} & 0 \\
        -h_k^{-2} & h_k^{-2} + h_{k+1}^{-2} & -h_{k+1}^{-2} \\
        0 & -h_{k+1}^{-2} & h_{k+1}^{-2}
	\end{bmatrix}\begin{bmatrix}
	    u_{k-1} \\ u_k \\ u_{k+1}
	\end{bmatrix}.
\end{equation}

\renewcommand{\KEeps}{\mathbf K_{\eps}^E}

\subsubsection{Topological derivative for the \enrichedMethod}\label{sec::discreteTopEnriched}
In this section, we compute the topological derivative for the \enrichedMethod. The results from \cite{gangl2023unified} used in \Cref{sec::discreteTopStandard} cannot be used here, because the linear system gets enlarged by the topological perturbation. In the following we use a direct approach (material derivative) rather than a Lagrangian approach. The following holds for interior nodes $x_k$. On the boundary the derivative of the enrichment function is unbounded and therefore needs a different treatment. For later use we prove the following matrix L'Hôspital rule.

\begin{lemma}\label{lemma::LHospital}
   Let $m,n \in \mathbb N$ and $\overline{\varepsilon}>0$. Let $A:[0,\overline{\varepsilon}) \rightarrow \mathbb R^{n\times n}$, $b:[0,\overline{\varepsilon}) \rightarrow \mathbb R^{n\times m}$ such that $A(0) = 0$, $b(0)=0$ and such that $A'_+(0):=\underset{\varepsilon \searrow 0}{\mbox{lim }} A(\varepsilon)/ \varepsilon$ and $b'_+(0):=\underset{\varepsilon \searrow 0}{\mbox{lim }} b(\varepsilon) / \varepsilon$ exist and $A'_+(0)$ is invertible.
Then $A(\varepsilon)$ is invertible for all $\varepsilon>0$ sufficiently small and
\begin{align}
    \underset{\varepsilon \searrow 0}{\mbox{lim }} A(\varepsilon)^{-1} b(\varepsilon) = A'_+(0) ^{-1}  b'_+(0) .
\end{align}
\end{lemma}
\begin{proof}
    Due to the right differentiability of $A$ and $b$ and $A(0)=b(0)=0$, we have
\begin{align*}
    A(\varepsilon) = \varepsilon ( A'_+(0) + r_A(\varepsilon) ), \quad b(\varepsilon) = \varepsilon ( b'_+(0) +  r_b(\varepsilon) ),
\end{align*}
with $\mathbb R^{n\times n} \ni r_A(\varepsilon)\rightarrow 0$ and $\mathbb R^{n\times m} \ni r_b(\varepsilon) \rightarrow 0$ as $\varepsilon \searrow 0$. Since, by assumption, $A'_+(0)$ is invertible and the set of invertible matrices is open, also $A'_+(0) + r_A(\varepsilon) $ is invertible for all sufficiently small $\varepsilon>0$. Thus, we have
\begin{align*}
    A(\varepsilon)^{-1} b(\varepsilon) = \left( A'_+(0) + r_A(\varepsilon) \right)^{-1} \left( b'_+(0) + r_b(\varepsilon) \right) \rightarrow A'_+(0)^{-1} b'_+(0)
\end{align*}
by the continuity of matrix inversion.
\end{proof}

For $x \in D$, let $ \mathbf N^S(x) = [N^S_1(x), \dots,N^S_{n_S}(x)  ]^\top \in \mathbb R^{n_S}$ and $ \mathbf N^E_\eps(x) = [N^E_{\eps,1}(x), \dots,N^E_{\eps,n_E}(x)  ]^\top \in \mathbb R^{n_E}$ the vectors of standard and enriched basis functions, respectively, such that the unperturbed state $u$ is represented by
\begin{equation}
	u^h(x) = \mathbf N^S(x)^\top \mathbf u = \sum_{i=1}^{n_S} N^S_i(x) u_i,
\end{equation}
and the perturbed state $u^h_\eps$ is represented by
\begin{equation} \label{eq_ueps_x}
	u^h_\eps(x) = \mathbf N^S(x)^\top \mathbf u^S_\eps + \mathbf N^E_\eps(x)^\top \mathbf u^E_\eps.
\end{equation}
Recall that the vector $\mathbf u\in \mathbb R^{n_S}$ is the solution of
\begin{equation} \label{eq_Ksu_f}
    \mathbf K^S \mathbf u = \mathbf f^S,
\end{equation}
whereas the vectors $\mathbf u^S_\eps\in \mathbb R^{n_S}$ and  $\mathbf u^E_\eps\in \mathbb R^{n_E}$ are the solution of
\begin{equation}\label{eq::topDirLinearSystem}
    \begin{bmatrix}
			\KmatSeps & \KmatSEeps \\
			\KmatESeps &  \KEeps
		\end{bmatrix}
  \begin{bmatrix}
      \mathbf u^S_\eps \\ \mathbf u^E_\eps
  \end{bmatrix} =
  \begin{bmatrix}
      \mathbf f^S \\ \mathbf f^E_\eps
  \end{bmatrix}.
\end{equation}
Note that, in the present case of a topological perturbation, for $\eps\searrow0$ the system \eqref{eq::topDirLinearSystem}
becomes singular. In particular the $n_E=2$ enrichment functions become identical and they are in $\text{span}\{N_i^S: i=1,\dots, n_S \}$, i.e., they can be represented by the standard basis functions. Thus, there exists a matrix $\mathbf H \in \R^{n_E\times n_S}$ such that
\begin{equation} \label{eq_NE_HNS}
    \mathbf N^E_0(x) = \mathbf H \mathbf N^S(x),
\end{equation}
for all $x \in D$ and it follows
\begin{align}\label{eq::HFFb}
    \mathbf H \mathbf K^S_0 &= \mathbf K^{ES}_0, \quad    \mathbf H\mathbf K^{SE}_0 = \mathbf K^{E}_0, \quad
    \mathbf H \mathbf f^S = \mathbf f^E_0.
\end{align}
In the considered case, i.e. the point of perturbation is a node $x_k$ and the polynomial degree is $p=1$, the matrix $\mathbf H$ has the entries $\mathbf H[i,j] =N_{0,i}^E(x_j)$. Furthermore,
note that $\lim_{\eps\searrow 0}u_\eps^h = u^h$, but it will turn out that $\lim_{\eps\searrow 0}\mathbf u_\eps^E \neq \mathbf 0$ and thus $\lim_{\eps\searrow 0}\mathbf u^S_\eps \neq \mathbf u$.
This observation motivates to introduce $\mathbf u_\eps \in \R^{n_S}$ as the solution of the perturbed system arising from the standard discretization,
\begin{align}\label{eq::perturbedUs}
	\mathbf K^S_\eps \mathbf u_\eps &= \mathbf f^S.
\end{align}
From the first line of the system \eqref{eq::topDirLinearSystem} and \eqref{eq::perturbedUs} we obtain
\begin{equation*}
    \mathbf u_\eps^S - \mathbf u_\eps = - (\KmatSeps)^{-1} \KmatSEeps \mathbf u^E_\eps = -\mathbf Q_\eps \mathbf u_\eps^E
\end{equation*}
with $\mathbf Q_\eps :=  (\KmatSeps)^{-1} \KmatSEeps \in \R^{n_S\times n_E}$.
Therefore, we have the decomposition
\begin{equation}\label{eq::decomposition}
    \mathbf u^S_\eps = \mathbf u_\eps-\mathbf Q_\eps \mathbf u_\eps^E.
\end{equation}

\begin{lemma} Assume that $\mathbf G:={d\mathbf K^E} - {d\mathbf K^{ES}}\mathbf H^\top + \mathbf H {d\mathbf K^S}\mathbf H^\top - \mathbf H d\mathbf K^{SE}$ is invertible. Then, the vector $\mathbf u^E_0 := \lim_{\eps\searrow0} \mathbf u_\eps^E$ can be computed by
\begin{equation}
	\mathbf u^E_0 =  \mathbf G^{-1}\left(d\mathbf f^E -\left(d\mathbf K^{ES} - \mathbf H d\mathbf K^S\right)\mathbf u \right).
\end{equation}
\end{lemma}
\begin{proof}
By means of the Schur complement of \eqref{eq::topDirLinearSystem} we have
\begin{align}\label{eq::SchurComplement}
	(\mathbf K^E_\eps - \KmatESeps(\KmatSeps)^{-1} \KmatSEeps)  \mathbf u^E_\eps &= \mathbf f^E_\eps - \KmatESeps (\KmatSeps)^{-1} \mathbf f^S.
\end{align}
Thus, the we will apply Lemma \ref{lemma::LHospital} with $\mathbf A(\eps) = \mathbf K^E_\eps - \KmatESeps(\KmatSeps)^{-1}\KmatSEeps$ and $\mathbf b(\eps)= \mathbf f^E_\eps -  \KmatESeps(\KmatSeps)^{-1}\mathbf f^S$ to derive $\uu_0^E$.
Let $\overline \eps=\min(x_k - x_{k-1},x_{k+1} - x_{k})$ be the local mesh width for node $x_k$. For $0<\eps<\overline \eps$ the basis functions $\psi_i$ are linearly independent, and therefore the system matrix in \eqref{eq::topDirLinearSystem} as well as the matrix $\KmatSeps$ are invertible. It follows that the Schur complement matrix $\mathbf A_\eps$ is invertible on $(0,\overline \eps)$.
 Note that by \eqref{eq::HFFb} we have
    \begin{equation}
        \lim_{\eps \rightarrow 0}\mathbf A_\eps =  \mathbf K^E_0 - \KK^{ES}_0(\KK^S_0)^{-1}\KK^{SE}_0= K^E_0 - \mathbf H \KK^{SE}_0 = \mathbf 0,
    \end{equation}
    and
    \begin{equation}
    \lim_{\eps \rightarrow 0} \mathbf b_\eps = \mathbf f^E_0 -  \KK^{ES}_0(\KK^S_0)^{-1}\mathbf f^S = \mathbf f^E_0 -  \mathbf H\mathbf f^S = \mathbf 0.
    \end{equation}
    Moreover, the mappings $\eps \mapsto \mathbf K^{ES}_\eps$, $\eps \mapsto \mathbf K^{S}_\eps$, $\eps \mapsto \mathbf K^{SE}_\eps$, $\eps \mapsto \mathbf K^{E}_\eps$, $\eps \mapsto \mathbf f^{E}_\eps$ are continuously differentiable on $[0, \overline \eps)$.
    We illustrate this for the first column of $\mathbf K^{SE}_\eps$: For $\eps > 0$ and $\eps = 0$, it holds \begin{align} \nonumber
        \KK_\eps^{SE}[i,1] =& \int_0^{\kappa-\eps} \lambda_1 \frac{d}{dx}N^S_i  \frac{d}{dx}N^E_{1, \eps}\ddx + \int_{\kappa-\eps}^{\kappa+\eps} \lambda_2 \frac{d}{dx}N^S_i  \frac{d}{dx}N^E_{1, \eps} + \int_{\kappa+\eps}^\ell \lambda_1 \frac{d}{dx}N^S_i  \frac{d}{dx}N^E_{1, \eps}\ddx\\
        &=\int_0^{\kappa-\eps} \lambda_1 \frac{d}{dx}N^S_i  \frac{1}{\kappa - \eps}dx +\int_{\kappa-\eps}^{\kappa+\eps} \lambda_2 \frac{d}{dx}N^S_i  \frac{1}{\ell -(\kappa - \eps)}dx + \int_{\kappa+\eps}^\ell \lambda_1\frac{d}{dx}N^S_i  \frac{1}{\ell-(\kappa-\eps)}\ddx, \label{eq_K_eps_SE}\\
        \KK_0^{SE}[i,1] =&\int_0^{\kappa} \lambda_1 \frac{d}{dx}N^S_i  \frac{1}{\kappa}dx + \int_{\kappa}^\ell \lambda_1\frac{d}{dx}N^S_i  \frac{1}{\ell-\kappa}\ddx, \nonumber
    \end{align}
    respectively, and in particular, by splitting the integration domain into the subintervals $(0, \kappa-\eps)$, $(\kappa+\eps, \ell)$, $(\kappa-\eps, \kappa)$ and $(\kappa-\eps, \kappa+\eps)$ it can be seen that
    \begin{align*}
        \frac{1}{\eps}(\KK_\eps^{SE}[i,1] - \KK_0^{SE}[i,1]) \rightarrow& \frac{1}{\kappa^2}\int_0^\kappa \lambda_1 \frac{d}{dx} N_i^S \,dx - \frac{1}{(\ell-\kappa)^2}\int_\kappa^\ell \lambda_1 \frac{d}{dx} N_i^S \,dx \\
        &+ \lambda_2 \frac{d}{dx}N_i^S(\kappa^-) \frac{1}{\ell-\kappa} - \lambda_1 \frac{d}{dx}N_i^S(\kappa^-) \frac{1}{\kappa}\\
        &+ \lambda_2 \frac{d}{dx}N_i^S(\kappa^+) \frac{1}{\ell-\kappa} - \lambda_1 \frac{d}{dx}N_i^S(\kappa^+) \frac{1}{\ell-\kappa} =: d \KK^{SE}[i, 1]
    \end{align*}
    as $\eps \rightarrow 0$. On the other hand, by applying the Leibniz rule
    $$\frac{d}{d\eps} \int_{a(\eps)}^{b(\eps)} f(\eps, x) \ddx = f(\eps, b(\eps)) b'(\eps) - f(\eps, a(\eps))a'(\eps) + \int_{a(\eps)}^{b(\eps)} \frac{\partial}{\partial \eps} f(\eps, x) \ddx$$
    to the three integrals in \eqref{eq_K_eps_SE} we have
    \begin{align*}
        \frac{d}{d\eps}&\KK_\eps^{SE}[i,1] =
        -\lambda_1 \frac{d}{dx}N^S_i( \kappa-\eps)  \frac{1}{\kappa - \eps} + \int_0^{\kappa - \eps} \lambda_1 \frac{d}{dx}N^S_i  \frac{1}{(\kappa - \eps)^2} \ddx \\
        &+ \lambda_2 \frac{d}{dx}N^S_i(\kappa+\eps)  \frac{1}{\ell -(\kappa - \eps)} + \lambda_2 \frac{d}{dx}N^S_i(\kappa-\eps)  \frac{1}{\ell -(\kappa - \eps)} + \int_{\kappa-\eps}^{\kappa+\eps} \lambda_2 \frac{d}{dx}N^S_i  \frac{-1}{(\ell -(\kappa - \eps))^2} \ddx \\
        &-\lambda_1\frac{d}{dx}N^S_i(\kappa +\eps)  \frac{1}{\ell-(\kappa-\eps)} + \int_{\kappa + \eps}^\ell \lambda_1\frac{d}{dx}N^S_i  \frac{-1}{(\ell-(\kappa-\eps))^2} \ddx.
    \end{align*}
    Since the integral in the second line above vanishes due to the symmetry of $N_i^S$ around the node $\kappa$, it follows that $\lim_{\eps \rightarrow 0} \frac{d}{d\eps} \KK_\eps^{SE}[i,1] = d \KK^{SE}[i,1]$. Since $i \in \{1, \dots, n_S\}$ was arbitrary and the same can be conducted for the second column, we see that $\eps \mapsto \KK^{SE}_\eps$ is continuously differentiable at $\eps=0$. Similar arguments for the other mappings mentioned above yield that also the composed mappings $\eps \mapsto \mathbf A_\eps$ and $\eps \mapsto \mathbf b_\eps$ are continuously differentiable on $[0, \overline \eps)$.

    Thus, employing Lemma \ref{lemma::LHospital} yields
\begin{equation}\label{eq::uE0}
	\mathbf u^E_0 = \left(\lim_{\eps\searrow 0}\frac{d}{d\eps}\left(\mathbf K^E_\eps - \KmatESeps(\KmatSeps)^{-1}\KmatSEeps\right)\right)^{-1}\left(\lim_{\eps\searrow 0}\frac{d}{d\eps}\left(\mathbf f^E_\eps -\KmatESeps(\KmatSeps)^{-1}\mathbf f^S\right) \right),
\end{equation}
provided that these limits exist and that the first limit is invertible, which will be computed in the following.
Note that, due to $\KK_\eps^{ES} = (\KK_\eps^{SE})^\top$ and $\KK_\eps^S = (\KK_\eps^S)^\top$, it holds
\begin{equation}\label{eq::defQtilde}
    \Qtilde = \KmatESeps(\KmatSeps)^{-1}.
\end{equation}
Using
\begin{equation}
    \frac{d}{d\eps}(\mathbf A_{\eps}^{-1}) = -\mathbf A_{\eps}^{-1} \frac{d\mathbf A_{\eps}}{d\eps}\mathbf A_{\eps}^{-1},
\end{equation}
we can compute
\begin{align*}
    \frac{d}{d\eps}\left(\KmatESeps(\KmatSeps)^{-1}\KmatSEeps\right) &=\frac{d\KmatESeps}{d\eps}(\KmatSeps)^{-1}\KmatSEeps + \KmatESeps\frac{d(\KmatSeps)^{-1}}{d\eps}\KmatSEeps+ \KmatESeps(\KmatSeps)^{-1}\frac{d\KmatSEeps}{d\eps} \\
    &=\frac{d\KmatESeps}{d\eps}(\KmatSeps)^{-1}\KmatSEeps - \KmatESeps(\KmatSeps)^{-1}\frac{d\KmatSeps}{d\eps}(\KmatSeps)^{-1}\KmatSEeps+ \KmatESeps(\KmatSeps)^{-1}\frac{d\KmatSEeps}{d\eps}\\
    &=\frac{d\KmatESeps}{d\eps}\mathbf Q_\eps - \Qtilde\frac{d\KmatSeps}{d\eps}\mathbf Q_\eps+ \Qtilde\frac{d\KmatSEeps}{d\eps},
\end{align*}
and
\begin{align*}
    \frac{d}{d\eps}\left(\KmatESeps(\KmatSeps)^{-1}\right) &=\frac{d\KmatESeps}{d\eps}(\KmatSeps)^{-1} + \KmatESeps\frac{d(\KmatSeps)^{-1}}{d\eps} \\
    &=\frac{d\KmatESeps}{d\eps}(\KmatSeps)^{-1} - \KmatESeps(\KmatSeps)^{-1}\frac{d\KmatSeps}{d\eps}(\KmatSeps)^{-1}\\
    &=\left(\frac{d\KmatESeps}{d\eps} - \Qtilde\frac{d\KmatSeps}{d\eps}\right)(\KmatSeps)^{-1}.
\end{align*}
Note that from \eqref{eq::defQtilde} and the definition of $\Q$ it follows,
\begin{align}
\mathbf H = \mathbf K_0^{ES}(\mathbf K_0^S)^{-1} = \mbox{lim}_{\eps\searrow0}  \mathbf K_\eps^{ES}(\mathbf K_\eps^S)^{-1} = \mbox{lim}_{\eps\searrow0} \Qtilde,\\
\mathbf H^\top = (\mathbf K_0^S)^{-1}\mathbf K_0^{SE} = \mbox{lim}_{\eps\searrow0}  (\mathbf K_\eps^S)^{-1}\mathbf K_\eps^{SE} = \mbox{lim}_{\eps\searrow0} \Q. \label{eq_Qeps_Ht}
\end{align}
Using that $(\KK^S_\eps)^{-1}\ff^S = \uu_\eps \rightarrow \uu$ defined in \eqref{eq_Ksu_f} yields
\begin{align*}
    \mathbf u^E_0 = \left({d\mathbf K^E} - {d\mathbf K^{ES}}\mathbf H^\top + \mathbf H {d\mathbf K^S}\mathbf H^\top - \mathbf H d\mathbf K^{SE}\right)^{-1}\left(d\mathbf f^E -\left(d\mathbf K^{ES} - \mathbf H d\mathbf K^S\right)\mathbf u \right)
\end{align*}
where we have used the continuous differentiability of all sub-matrices and sub-vectors of the system \eqref{eq::topDirLinearSystem} with respect to $\eps$ at $\eps=0$. Since the matrix is invertible by assumption, this completes the proof.

\end{proof}

\begin{remark}
    We remark that we have checked that the $2\times 2$ matrix $d\mathbf K^E - {d\mathbf K^{ES}}\mathbf H^\top + \mathbf H {d\mathbf K^S}\mathbf H^\top - \mathbf H d\mathbf K^{SE}$ is invertible by means of symbolic computations done in Matlab for the special cases $n_S=2,4,8$.
\end{remark}

 For the next result we recall
\begin{align*}
    d\mathbf M^{SE} &= \int_0^\ell \mathbf N^S \left(\left. \frac{d\mathbf N^E_\eps}{d\eps}\right\rvert_{\eps=0} \right)^\top \ddx \in \R^{n_S\times n_E},\\
    d\mathbf m^E &= \int_0^\ell \hat u(x) \left. \frac{d\mathbf N^E_\eps}{d\eps}\right\rvert_{\eps=0} \ddx \in \R^{n_E}.
\end{align*}
\begin{theorem}
The topological derivative for the \enrichedMethod at node $x_k$ is given by
\begin{equation}\label{eq::topological_derivativeC}
	d_T\bar G(x_k) = S_0 + S_1  + S_2,
\end{equation}
with
\begin{subequations}
	\begin{align}
		S_0 &=  \frac12 \mathbf p^\top d\mathbf K^S \mathbf u,
        \\
		S_1 &=  (\mathbf u^\top d\mathbf M^{SE}-(d\mathbf m^E )^\top)\mathbf u^E_0, \\
		S_2 &= \frac12 \mathbf p^\top \left(d\mathbf K^S\mathbf H^\top - d\mathbf K^{SE} \right) \,\mathbf u^E_0,
	\end{align}
\end{subequations}
where the adjoint state $\mathbf p\in \R^{n_S}$ is the solution of
\begin{align}
	(\mathbf K^S_0)^\top \mathbf p &= -2(\mathbf M \mathbf u - \mathbf m) 
 .
 \end{align}
\end{theorem}
\begin{proof}
In 1D we can interpret the topological derivative as ordinary derivative,
\begin{equation*}
	d_T\bar G(x_k)=\lim_{\eps\searrow 0} \frac{G(u^h_\eps)-G(u^h)}{2\eps} = \frac{1}{2}\lim_{\eps\searrow 0}\frac{dG(u^h_\eps)}{d\eps} = \int_0^\ell (u^h-\hat u(x)) \left. \frac{du^h_\eps}{d\eps}\right \rvert_{\eps=0} \ddx.
\end{equation*}
Using \eqref{eq_ueps_x}, the derivative of the state can be written as
\begin{equation}
\begin{aligned}
	\left. \frac{du_\eps^h}{d\eps}\right \rvert_{\eps=0} &= \lim_{\eps\searrow 0} \frac{(\mathbf N^S)^\top \mathbf u^S_\eps + (\mathbf N^E_\eps)^\top \mathbf u^E_\eps - (\mathbf N^S)^\top \mathbf u}{\eps} \\
	&=  (\mathbf N^S)^\top \lim_{\eps\searrow 0}\frac{\mathbf u_\eps - \mathbf u}{\eps} + \lim_{\eps\searrow 0} \frac{\left((\mathbf N^E_\eps)^\top-(\mathbf N^S)^\top\mathbf Q_\eps\right)}{\eps}\mathbf u^E_\eps, \label{eq::derivativeU}
\end{aligned}
\end{equation}
where we have used the decomposition
\eqref{eq::decomposition}.
Adding and subtracting $\tfrac{1}{\eps} (\mathbf N^E_0)^\top\mathbf u^E_0$ to \eqref{eq::derivativeU} allows to reformulate it to
\begin{align*}
	\left. \frac{du^h_\eps}{d\eps}\right \rvert_{\eps=0} &= (\mathbf N^S)^\top \lim_{\eps\searrow 0}\frac{\mathbf u_\eps - \mathbf u}{\eps}
	+
	\lim_{\eps\searrow 0} \frac{(\mathbf N^E_\eps)^\top-(\mathbf N^E_0)^\top}{\eps}\mathbf u^E_0
	+
	\lim_{\eps\searrow 0} \frac{\left((\mathbf N^E_0)^\top-(\mathbf N^S)^\top\mathbf Q_\eps\right)}{\eps}\mathbf u^E_0 \\
	&= (\mathbf N^S)^\top \frac{d\mathbf u_\eps}{d\eps}\bigg  \rvert_{\eps=0} + \left(\frac{d\mathbf N^E_\eps}{d\eps}^\top\bigg  \rvert_{\eps=0}- (\mathbf N^S)^\top\frac{d\mathbf Q_\eps}{d\eps}\bigg \rvert_{\eps=0}\right)\mathbf u^E_0
\end{align*}
where we also used \eqref{eq_Qeps_Ht} and \eqref{eq_NE_HNS}.
Subtracting the unperturbed state equation from \eqref{eq::perturbedUs} gives us
\begin{align*}
	\mathbf K^S_\eps \mathbf u_\eps - \mathbf K^S_0 \mathbf u &= (\mathbf K^S_\eps-\mathbf K^S_0) \mathbf u_\eps + \mathbf K^S_0 (\mathbf u_\eps-\mathbf u) = \mathbf 0,
\end{align*}
and we conclude that $\frac{d\mathbf u_\eps}{d\eps} =\lim_{\eps\searrow 0}\frac{\mathbf u_\eps-\mathbf u}{\eps}$ is the solution of
\begin{align}\label{eq::ueps}
	\mathbf K^S_0 \frac{d\mathbf u_\eps}{d\eps} = -\lim_{\eps\searrow 0}\frac{\mathbf K^S_\eps-\mathbf K^S_0}{\eps} \mathbf u = - d\mathbf K^S \mathbf u.
\end{align}
Furthermore,
\begin{align*}
    \frac{d\mathbf Q_\eps}{d\eps} &= \frac{d(\KmatSeps)^{-1}}{d\eps} \KmatSEeps + (\KmatSeps)^{-1} \frac{d\KmatSEeps}{d\eps}
    \\
    &=-(\KmatSeps)^{-1} \frac{d\KmatSeps}{d\eps} (\KmatSeps)^{-1}\KmatSEeps + (\KmatSeps)^{-1} \frac{d\KmatSEeps}{d\eps} \\
    &= (\KmatSeps)^{-1}\left(-\frac{d\KmatSeps}{d\eps} \Q+\frac{d\KmatSEeps}{d\eps}\right).
\end{align*}
Therefore, we can write for the topological derivative
\begin{align*}\label{eq::topological_derivativeA}
	d_T\bar G(x_k)&= \int_0^\ell (u^h-\hat u(x)) \left( (\mathbf N^S)^\top \frac{d\mathbf u_\eps}{d\eps}\bigg \rvert_{\eps=0} + \left(\frac{d\mathbf N^E_\eps}{d\eps}^\top\bigg \rvert_{\eps=0}- (\mathbf N^S)^\top\frac{d\mathbf Q_\eps}{d\eps}\bigg \rvert_{\eps=0}\right)\mathbf u^E_0 \right)\ddx \\
 &= \frac{1}{2}\int_0^\ell (-2)(u^h-\hat u(x))  (\mathbf N^S)^\top  \ddx \, (\mathbf K^S)^{-1}\; d\mathbf K^S \mathbf{u} \\
 & + \int_0^\ell (u^h-\hat u(x))  \frac{d\mathbf N^E_\eps}{d\eps}^\top\bigg \rvert_{\eps=0}\ddx \, \mathbf{u}_0^E\\
 & - \frac{1}{2}\int_0^\ell (-2)(u^h-\hat u(x))  (\mathbf N^S)^\top  \ddx \, (\mathbf K^S)^{-1}\; \left(-d\mathbf K^S\mathbf H^\top + d\mathbf K^{SE} \right) \mathbf{u}_0^E
\end{align*}
Using the adjoint state finishes the proof.
\end{proof}

\newcommand{\deLepsA}{\frac{1}{x_k-\eps}}
\newcommand{\deRepsA}{\frac{1}{\ell-x_k+\eps}}
\newcommand{\deLepsB}{\frac{1}{x_k+\eps}}
\newcommand{\deRepsB}{\frac{1}{\ell-x_k-\eps}}
\newcommand{\deLA}{\frac{1}{x_k}}
\newcommand{\deRA}{-\frac{1}{\ell-x_k}}
\newcommand{\deLB}{\frac{1}{x_k}}
\newcommand{\deRB}{-\frac{1}{\ell-x_k}}
\newcommand{\ddeLa}{\frac{1}{x_k^2}}
\newcommand{\ddeRa}{\frac{1}{(\ell-x_k)^2}}
\newcommand{\ddeLb}{-\frac{1}{x_k^2}}
\newcommand{\ddeRb}{-\frac{1}{(\ell-x_k)^2}}
\newcommand{\mat}[2]{\begin{bmatrix} #1 & #2 \end{bmatrix}}

\section{Impact of the discretization on smoothness of the shape derivative}
In this section, we present a theorem regarding the smoothness in space of the shape derivative. It turns out that the smoothness of the discretization (\ie the polynomial degree $p$ of the splines) determines the smoothness of the shape derivative.
\begin{theorem}\label{theorem::shape}
    Let $f\in C^{p-2}(D)$ and $\targetU \in C^{p-1}(D)$. Then, the mapping $\kappa \mapsto d_S\bar G(\kappa)$ is in $C^{p-2}(D)$.
\end{theorem}
\begin{proof}
Following a direct approach the shape derivative reads
\begin{align}
    d_S\bar G(\kappa) = \frac{d}{d\kappa}G(u^h(x,\kappa)) = 2 \int_0^\ell (u^h(x,\kappa)-\targetU(x)) \frac{d}{d\kappa}(u^h(x,\kappa)) \ddx
\end{align}
and thus the smoothness with respect to $\kappa$ of $\frac{d}{d\kappa}G(u^h(x,\kappa))$ is determined by the smoothness of $\frac{d}{d\kappa}u^h(x,\kappa)$. In particular, if $(x,\kappa)\mapsto u^h(x,\kappa) \in C^{p-1}(D)\times C^g(D)$ for some $g$ then $(x,\kappa)\mapsto (u^h(x,\kappa)-\targetU(x))\frac{d}{d\kappa}(u^h(x,\kappa)) \in C^{p-1}(D)\times C^{g-1}(D)$ and thus $\kappa \mapsto d_S\bar G(\kappa) \in C^{g-1}(D)$.

Recall that for the \standardMethod $u^h(x,\kappa) = \mathbf N^S(x) \mathbf u(\kappa)$, where $\mathbf u(\kappa)$ is given by
\begin{equation*}
    \mathbf u(\kappa) = \mathbf K^S(\kappa)^{-1} \mathbf f(\kappa).
\end{equation*}
We have $N_i \in C^{p-1}(D)$ and $\frac{dN_i}{dx} \in C^{p-2}(D)$.
The stiffness matrix is given by
\begin{align} \nonumber
    \mathbf K^S(\kappa)[i,j] &= \int_0^\ell\lambda_\Omega \frac{dN_i}{dx}(x) \frac{dN_j}{dx}(x) \ddx \\ &=\int_0^\kappa \lambda_1 \frac{dN_i}{dx}(x) \frac{dN_j}{dx}(x) \ddx + \int_\kappa^\ell \lambda_2 \frac{dN_i}{dx}(x) \frac{dN_j}{dx}(x) \ddx. \label{eq_KS_split}
\end{align}
Since $x \mapsto \frac{dN_i}{dx}(x) \frac{dN_j}{dx}(x) \in C^{p-2}(D)$, the primitive functions of the integrands in \eqref{eq_KS_split} are in $C^{p-1}$. Thus, by the fundamental theorem of calculus, we conclude that $\kappa \mapsto \mathbf K^S(\kappa)[i,j]\in C^{p-1}(D)$. Due to $f \in C^{p-2}(D)$, it holds $\kappa\mapsto\mathbf f(\kappa) \in C^{p-1}(D)$ and we obtain $g=p-1$.

For the \enrichedMethod we can use the explicit representation of the enrichment function \eqref{eq:enrichment} to get
\begin{align} \nonumber
    \KK^{SE}(\kappa)[i] &= \int_0^\kappa \lambda_1 \frac{dN_i}{dx}(x) \frac{1}{\kappa} \ddx + \int_\kappa^\ell \lambda_2 \frac{dN_i}{dx}(x) \frac{-1}{\ell-\kappa} \ddx \\
    &= \frac{1}{\kappa} \int_0^\kappa \lambda_1 \frac{dN_i}{dx}(x)  \ddx - \frac{1}{\ell-\kappa}\int_\kappa^\ell \lambda_2 \frac{dN_i}{dx}(x)  \ddx,
\end{align}
and
\begin{equation}
    K^E(\kappa) = \frac{\lambda_1}{\kappa} + \frac{\lambda_2}{\ell-\kappa}.
\end{equation}
Using the same argumentation as for the \standardMethod, we conclude $\kappa \mapsto \KK^{SE}(\kappa)[i] \in C^{p-1}(D)$. Moreover, it is readily seen that $\kappa \mapsto K^{E}(\kappa) \in C^\infty(D)$. Thus, the assertion follows.
\end{proof}

\section{Numerical results}
In this section, we present the results of numerical experiments and compare the continuous shape and topological derivatives recalled in \Cref{sec::continuous} with the discrete sensitivities obtained in \Cref{sec::numericalSensitivities}. The used parameters and functions are summarized in \Cref{tab::data}.
\begin{table}[htb]
	\centering
	\begin{tabular}{ccccc}
		\toprule
		$\ell$  &
		$\lambda_1$ &
		$\lambda_2$ &
		$f(x) $ &
		$\hat u(x)$ \\
		 $1$ & \corrections{$0.6$} & $0.2$ & $x$ & $(\ell-x)x$ \\
		\bottomrule
	\end{tabular}
	\caption{Setup for the numerical experiments}
	\label{tab::data}
\end{table}
\subsection{Discretization error for the direct state}
\begin{figure}[htb]
    \begin{subfigure}[t]{0.49\textwidth}
        \begin{tikzpicture}
            \begin{axis}[xmode=log, ymode=log,grid,
            xlabel={$m$}, ylabel={$L_2$-error},
            width=6.7cm,height=5.7cm,xmin=2,xmax=4096,]
        \addplot[dashed,domain=2:4096,black] {0.005/x)};    
        \addplot[col1,mark=o] table[col sep=comma,y index=1] {pic/error_convergence_FEM.txt};
        \addplot[col2,mark=+] table[col sep=comma,y index=2] {pic/error_convergence_FEM.txt};
        \addplot[col3,mark=square] table[col sep=comma,y index=3] {pic/error_convergence_FEM.txt};
        \legend{{$O(h)$}}
        \end{axis}
        \end{tikzpicture}
        \caption{$L_2$-error for \FEM}
        \label{fig:convergenceA}
    \end{subfigure}
     \begin{subfigure}[t]{0.49\textwidth}
        \begin{tikzpicture}
            \begin{axis}[xmode=log, ymode=log,grid,
            xlabel={$m$}, ylabel={$H^1$-error},
            width=6.5cm,height=5.7cm,xmin=2,xmax=4096,]

        \addplot[dashed,domain=2:4096,black] {0.1/(x^0.5))};
            
        \addplot[col1,mark=o] table[col sep=comma,y index=4] {pic/error_convergence_FEM.txt};
        \addplot[col2,mark=+] table[col sep=comma,y index=5] {pic/error_convergence_FEM.txt};
        \addplot[col3,mark=square] table[col sep=comma,y index=6] {pic/error_convergence_FEM.txt};
        
        \legend{{$O(h^{0.5})$}}
        \end{axis}
        \end{tikzpicture}
        \caption{$H^1$-error for \FEM}
        \label{fig:convergenceB}
    \end{subfigure}\vspace{0.5cm}
    \begin{subfigure}[t]{0.49\textwidth}
        \begin{tikzpicture}
            \begin{axis}[xmode=log, ymode=log,grid,
            xlabel={$m$}, ylabel={$L_2$-error},
            width=6.7cm,height=5.7cm,xmin=2,xmax=4096,]
        \addplot[dashed,domain=2:4096,black] {0.02/x^2)}; \addplot[domain=2:4096,black] {0.005/x^(2.5))};       
        \addplot[col1,mark=o] table[col sep=comma,y index=1] {pic/error_convergence_XFEM.txt};
        \addplot[col2,mark=+] table[col sep=comma,y index=2] {pic/error_convergence_XFEM.txt};
        \addplot[col3,mark=square] table[col sep=comma,y index=3] {pic/error_convergence_XFEM.txt};

        \legend{{$O(h^2)$},{$O(h^{2.5})$}}
        \end{axis}
        \end{tikzpicture}
        \caption{$L_2$-error for \XFEM}
        \label{fig:convergenceC}
    \end{subfigure}
     \begin{subfigure}[t]{0.49\textwidth}
        \begin{tikzpicture}
            \begin{axis}[xmode=log, ymode=log,grid,
            xlabel={$m$}, ylabel={$H^1$-error},
            width=6.5cm,height=5.7cm,xmin=2,xmax=4096,]

        \addplot[dashed,domain=2:4096,black] {0.1/(x^1))};
        \addplot[domain=2:4096,black] {0.01/(x^1.5))};
            
        \addplot[col1,mark=o] table[col sep=comma,y index=4] {pic/error_convergence_XFEM.txt};
        \addplot[col2,mark=+] table[col sep=comma,y index=5] {pic/error_convergence_XFEM.txt};
        \addplot[col3,mark=square] table[col sep=comma,y index=6] {pic/error_convergence_XFEM.txt};
        \legend{{$O(h)$},{$O(h^{1.5})$}}
        \end{axis}
        \end{tikzpicture}
        \caption{$H^1$-error for \XFEM}
        \label{fig:convergenceD}
    \end{subfigure}\vspace{0.4cm}
    \begin{subfigure}[t]{0.99\textwidth}
    \centering
    \begin{tikzpicture} 
    \begin{axis}[hide axis,
    legend style={draw=white!15!black,legend cell align=left,legend columns=-1,xmin=10,xmax=50,ymin=0,ymax=0.4,}    ]
    \addlegendimage{col1,mark=o}
    \addlegendentry{$p = 1$};
    \addlegendimage{col2,mark=+}
    \addlegendentry{$p = 2$};
    \addlegendimage{col3,mark=square}
    \addlegendentry{$p = 3$};
    \end{axis}
    \end{tikzpicture}
    \end{subfigure}
    \caption{Convergence of the $L_2$ and $H^1$ errors for different discretizations}
    \label{fig:convergenceStudy}
\end{figure}
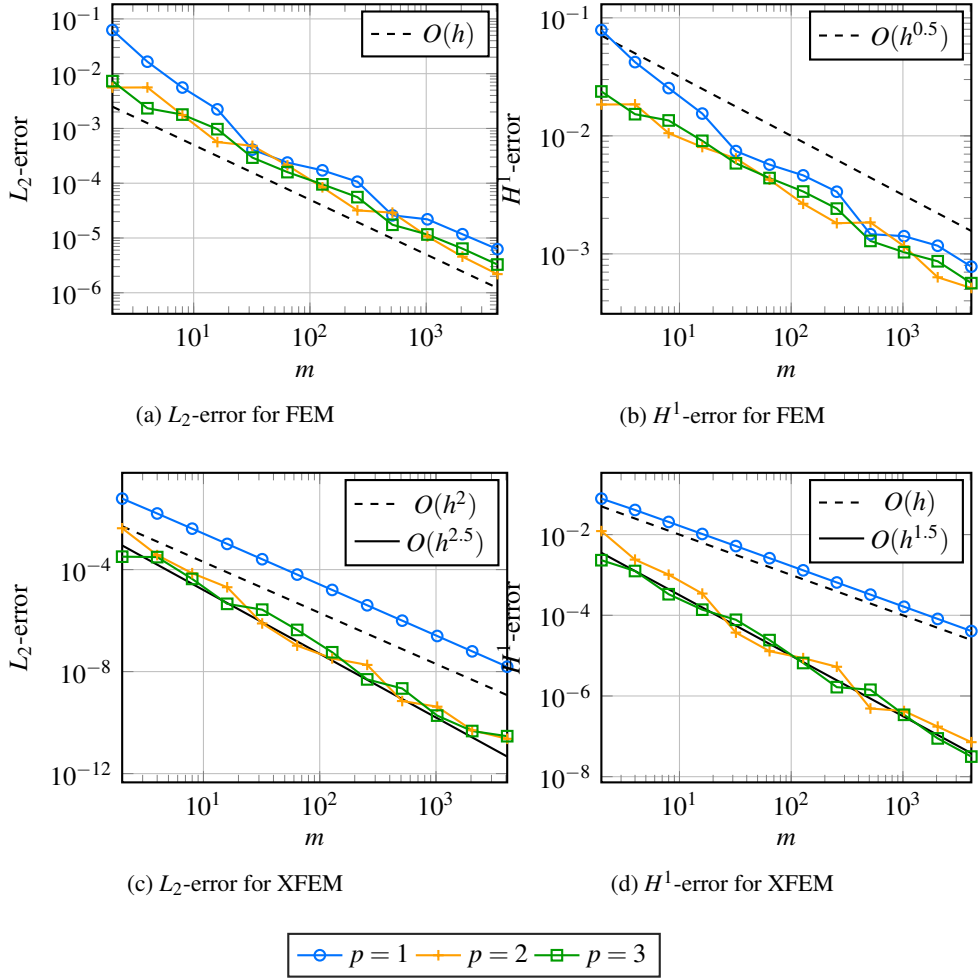

We begin by analyzing the convergence behavior of the $L_2$ and $H^1$ errors for the finite element approximation $u^h$ for the case of one interface at $\kappa = \sqrt{2}/5$, which is not resolved by the standard discretization. The results are presented in \Cref{fig:convergenceStudy}. For the \standardMethod, we observe the expected well known suboptimal convergence rates that are independent of the polynomial degree $p$. Specifically, the $L_2$ error exhibits a convergence rate of order $h$ (see \Cref{fig:convergenceA}), while the $H^1$ error converges at a rate of order $h^{0.5}$ (see \Cref{fig:convergenceB}).

In contrast, the \enrichedMethod results, shown in \Cref{fig:convergenceC} and \Cref{fig:convergenceD}, reveal a dependency on the polynomial degree. For the $L_2$ error (\Cref{fig:convergenceC}), we observe optimal convergence of order $h^2$ for $p=1$, and a higher rate of $h^{2.5}$ for $p > 1$. Similarly, for the $H^1$ error (\Cref{fig:convergenceD}), the convergence is optimal with order $h$ for $p=1$, while a rate of $h^{1.5}$ is obtained for $p > 1$. These suboptimal rates for higher polynomial degrees are attributed to the use of only linear enrichment functions in the \enrichedMethod formulation. Again, we remark that the choice of linear enrichment functions was made in order to facilitate the explicit computation of numerical sensitivities which allows to compare their behavior in different situations as presented in the following subsections.
\subsection{Influence of the smoothness of the discretization on the shape derivative}\label{sec:num_smoothness}
\begin{figure}[htb]
\begin{subfigure}[t]{0.49\textwidth}
    \begin{tikzpicture}
    \begin{axis}[grid, xlabel={$\kappa$}, ylabel={$G$},
        width=6.7cm,height=5.7cm,xmin=0,xmax=1,]
           
    \addplot[black] table[col sep=comma,y index=1] {pic/shape_analytic_solution_3.txt};
    \addplot[col1] table[col sep=comma,y index=1] {pic/shape_FEM_1_n4_3.txt};
    \addplot[col2] table[col sep=comma,y index=1] {pic/shape_FEM_2_n4_3.txt};
    \addplot[col3] table[col sep=comma,y index=1] {pic/shape_FEM_3_n4_3.txt};
    \end{axis}
    \end{tikzpicture}
    \caption{$m=4$}
    \label{fig:objectiveFEMA}
\end{subfigure}
\begin{subfigure}[t]{0.49\textwidth}
    \begin{tikzpicture}
    \begin{axis}[grid, xlabel={$\kappa$}, ylabel={$G$},
        width=6.7cm,height=5.7cm,xmin=0,xmax=1,legend pos=north west]
           
    \addplot[black] table[col sep=comma,y index=1] {pic/shape_analytic_solution_3.txt};
    \addplot[col1] table[col sep=comma,y index=1] {pic/shape_FEM_1_n8_3.txt};
    \addplot[col2] table[col sep=comma,y index=1] {pic/shape_FEM_2_n8_3.txt};
    \addplot[col3] table[col sep=comma,y index=1] {pic/shape_FEM_3_n8_3.txt};
    \legend{{analytic},$p=1$,$p=2$,$p=3$}
    \end{axis}
    \end{tikzpicture}
    \caption{$m=8$}
    \label{fig:objectiveFEMB}
\end{subfigure}
\caption{The objective function in dependency of the interface position for different discretizations using the \standardMethod}
\label{fig:objectiveFEM}
\end{figure}
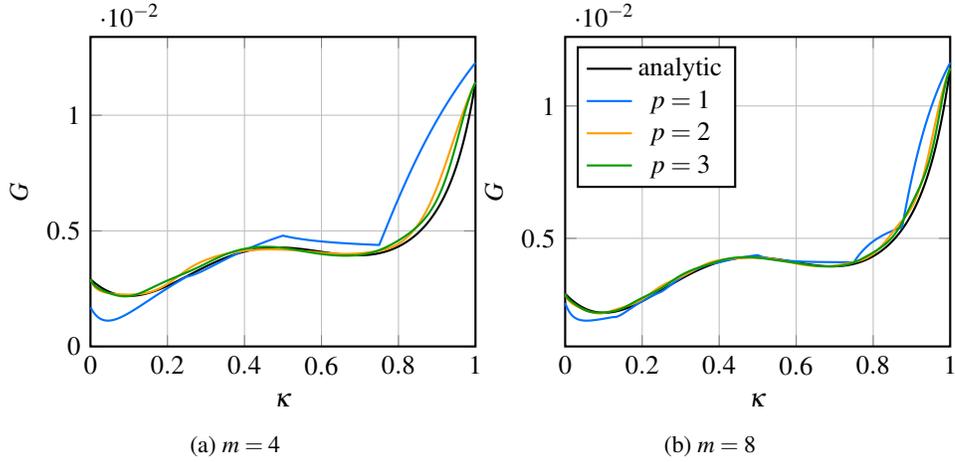
In this section, we numerically investigate the influence of the smoothness of the discretization on the shape derivative. To this end, we consider a one-interface configuration with the interface located at a varying position $\kappa$.

Before analyzing the shape derivative, we first examine the behavior of the objective function. The resulting objective values for various interface positions and discretization orders are compared to the analytical reference in \Cref{fig:objectiveFEM}. For $p = 1$, we observe noticeable kinks in the graph of the objective function at element boundaries, indicating low regularity. In contrast, for $p > 1$, the curves appear smoother.

As shown in \Cref{fig:objectiveFEMB}, where a mesh with 8 elements is used, the discrete solutions qualitatively match the analytical solution for $p > 1$. In \Cref{fig:shapeFEM}, the shape derivative \eqref{eq::derivativeTrackingShapeStandard} is plotted as a function of the interface position. For $p = 1$, jumps are clearly visible, whereas for $p > 1$, smoother curves are obtained, in agreement with \Cref{theorem::shape}.

In contrast to the objective function plotted in \Cref{fig:objectiveFEM}, however, the shape derivative exhibits oscillations even for $p > 1$, and these oscillations become more pronounced with mesh refinement.
\begin{figure}[htb]
\begin{subfigure}[t]{0.49\textwidth}
    \begin{tikzpicture}
    \begin{axis}[grid, xlabel={$\kappa$}, ylabel={$dG$},
        width=6.7cm,height=5.7cm,xmin=0,xmax=1,unbounded coords=jump,
        tick label style={/pgf/number format/fixed},]
           
    \addplot[black] table[col sep=comma,y index=2] {pic/shape_analytic_solution_3.txt};
    \addplot[col1] table[col sep=comma,y index=2] {pic/shape_FEM_1_n4_3.txt};
    \addplot[col2] table[col sep=comma,y index=2] {pic/shape_FEM_2_n4_3.txt};
    \addplot[col3] table[col sep=comma,y index=2] {pic/shape_FEM_3_n4_3.txt};
    \end{axis}
    \end{tikzpicture}
    \caption{$m=4$}
    \label{fig:shapeFEMA}
\end{subfigure}
\begin{subfigure}[t]{0.49\textwidth}
    \begin{tikzpicture}
    \begin{axis}[grid, xlabel={$\kappa$}, ylabel={$dG$},
        width=6.7cm,height=5.7cm,xmin=0,xmax=1, legend pos=north west,unbounded coords=jump,
        tick label style={/pgf/number format/fixed}, ]
           
    \addplot[black] table[col sep=comma,y index=2] {pic/shape_analytic_solution_3.txt};
    \addplot[col1] table[col sep=comma,y index=2] {pic/shape_FEM_1_n8_3.txt};
    \addplot[col2] table[col sep=comma,y index=2] {pic/shape_FEM_2_n8_3.txt};
    \addplot[col3] table[col sep=comma,y index=2] {pic/shape_FEM_3_n8_3.txt};
    \legend{{analytic},$p=1$,$p=2$,$p=3$}
    \end{axis}
    \end{tikzpicture}
    \caption{$m=8$}
    \label{fig:shapeFEMB}
\end{subfigure}
    \caption{The shape derivative in dependency of the interface position for different discretizations using the \standardMethod}
    \label{fig:shapeFEM}
\end{figure}

\subsection{Influence of the interface resolution on the shape derivative}
\begin{figure}[htb]
\begin{subfigure}[t]{0.49\textwidth}
    \begin{tikzpicture}
    \begin{axis}[grid, xlabel={$\kappa$}, ylabel={$dG$},
        width=6.7cm,height=5.7cm,xmin=0,xmax=1,unbounded coords=jump,
        tick label style={/pgf/number format/fixed},]
           
    \addplot[black] table[col sep=comma,y index=2] {pic/shape_analytic_solution_3.txt};
    \addplot[col1] table[col sep=comma,y index=2] {pic/shape_XFEM_1_n4_3.txt};
    \addplot[col2] table[col sep=comma,y index=2] {pic/shape_XFEM_2_n4_3.txt};
    \addplot[col3] table[col sep=comma,y index=2] {pic/shape_XFEM_3_n4_3.txt};
    \end{axis}
    \end{tikzpicture}
    \caption{$m=4$}
    \label{fig:shapeXFEMA}
\end{subfigure}
\begin{subfigure}[t]{0.49\textwidth}
    \begin{tikzpicture}
    \begin{axis}[grid, xlabel={$\kappa$}, ylabel={$dG$},
        width=6.7cm,height=5.7cm,xmin=0,xmax=1, legend pos=north west,unbounded coords=jump,
        tick label style={/pgf/number format/fixed},]
           
    \addplot[black] table[col sep=comma,y index=2] {pic/shape_analytic_solution_3.txt};
    \addplot[col1] table[col sep=comma,y index=2] {pic/shape_XFEM_1_n8_3.txt};
    \addplot[col2] table[col sep=comma,y index=2] {pic/shape_XFEM_2_n8_3.txt};
    \addplot[col3] table[col sep=comma,y index=2] {pic/shape_XFEM_3_n8_3.txt};
    \legend{{analytic},$p=1$,$p=2$,$p=3$}
    \end{axis}
    \end{tikzpicture}
    \caption{$m=8$}
    \label{fig:shapeXFEMB}
\end{subfigure}
    \caption{The shape derivative in dependency of the interface position for different discretizations based on the \enrichedMethod}
    \label{fig:shapeXFEM}
\end{figure}
\corrections{In section \ref{sec:num_smoothness}, the material interface was not resolved by the discretization.}
In this section, we comment on the influence of \corrections{resolving the interface by employing the \enrichedMethod} on the shape derivative. The results for different discretizations are shown in \Cref{fig:shapeXFEM}. Similar to the discretization used in the \standardMethod, we observe jumps in the shape derivative for $p = 1$. However, in this case, the jumps between elements diminish with mesh refinement. For $p > 1$, the shape derivative values obtained using the \enrichedMethod show good agreement with the analytical shape derivative obtained by differentiating \eqref{eq::G_kappa_explicit}.
\subsection{Comparison of continuous and numerical sensitivity formulas}
\begin{figure}[htb]
\begin{subfigure}[t]{0.49\textwidth}
    \begin{tikzpicture}
    \begin{axis}[grid, xlabel={$\kappa$}, ylabel={$dG$},
        width=6.7cm,height=5.7cm,xmin=0,xmax=1,unbounded coords=jump,
        tick label style={/pgf/number format/fixed},]

    \addplot[black] table[col sep=comma,y index=2] {pic/shape_analytic_solution_3.txt};
    \addplot[col1] table[col sep=comma,y index=2] {pic/shape_FEM_1_n4_3.txt};
    \addplot[col2] table[col sep=comma,y index=3] {pic/shape_FEM_2_n4_3.txt};
    \end{axis}
    \end{tikzpicture}
    \caption{\FEM, $m=4$, $p=1$}
    \label{fig:shapeFormulaA}
\end{subfigure}
\begin{subfigure}[t]{0.49\textwidth}
    \begin{tikzpicture}
    \begin{axis}[grid, xlabel={$\kappa$}, ylabel={$dG$},
        width=6.7cm,height=5.7cm,xmin=0,xmax=1,unbounded coords=jump,
        tick label style={/pgf/number format/fixed},]

    \addplot[black] table[col sep=comma,y index=2] {pic/shape_analytic_solution_3.txt};
    \addplot[col1] table[col sep=comma,y index=2] {pic/shape_FEM_1_n16_3.txt};
    \addplot[col2] table[col sep=comma,y index=3] {pic/shape_FEM_2_n16_3.txt};
    \end{axis}
    \end{tikzpicture}
    \caption{\FEM, $m=16$, $p=1$}
    \label{fig:shapeFormulaA2}
\end{subfigure}
\vspace{0.6cm}

\begin{subfigure}[t]{0.49\textwidth}
    \begin{tikzpicture}
    \begin{axis}[grid, xlabel={$\kappa$}, ylabel={$dG$},
        width=6.7cm,height=5.7cm,xmin=0,xmax=1, legend pos=north west,unbounded coords=jump,
        tick label style={/pgf/number format/fixed},]

    \addplot[black] table[col sep=comma,y index=2] {pic/shape_analytic_solution_3.txt};
    \addplot[col1] table[col sep=comma,y index=2] {pic/shape_XFEM_1_n4_3.txt};
    \addplot[col2] table[col sep=comma,y index=3] {pic/shape_XFEM_1_n4_3.txt};

    \end{axis}
    \end{tikzpicture}
    \caption{\XFEM, $m=4$, $p=1$}
    \label{fig:shapeFormulaB}
    \end{subfigure}
    \begin{subfigure}[t]{0.49\textwidth}
    \begin{tikzpicture}
    \begin{axis}[grid, xlabel={$\kappa$}, ylabel={$dG$},
        width=6.7cm,height=5.7cm,xmin=0,xmax=1, legend pos=north west,unbounded coords=jump,
        tick label style={/pgf/number format/fixed},]

    \addplot[black] table[col sep=comma,y index=2] {pic/shape_analytic_solution_3.txt};
    \addplot[col1] table[col sep=comma,y index=2] {pic/shape_XFEM_1_n16_3.txt};
    \addplot[col2] table[col sep=comma,y index=3] {pic/shape_XFEM_1_n16_3.txt};
    \legend{{analytic},{DP},{CP}}
    \end{axis}
    \end{tikzpicture}
    \caption{\XFEM, $m=16$, $p=1$}
    \label{fig:shapeFormulaB2}
    \end{subfigure}
    \caption{Shape derivative versus interface position $\kappa$ for different sensitivity formulas and discretizations.}
    \label{fig:shapeFormula}
\end{figure}
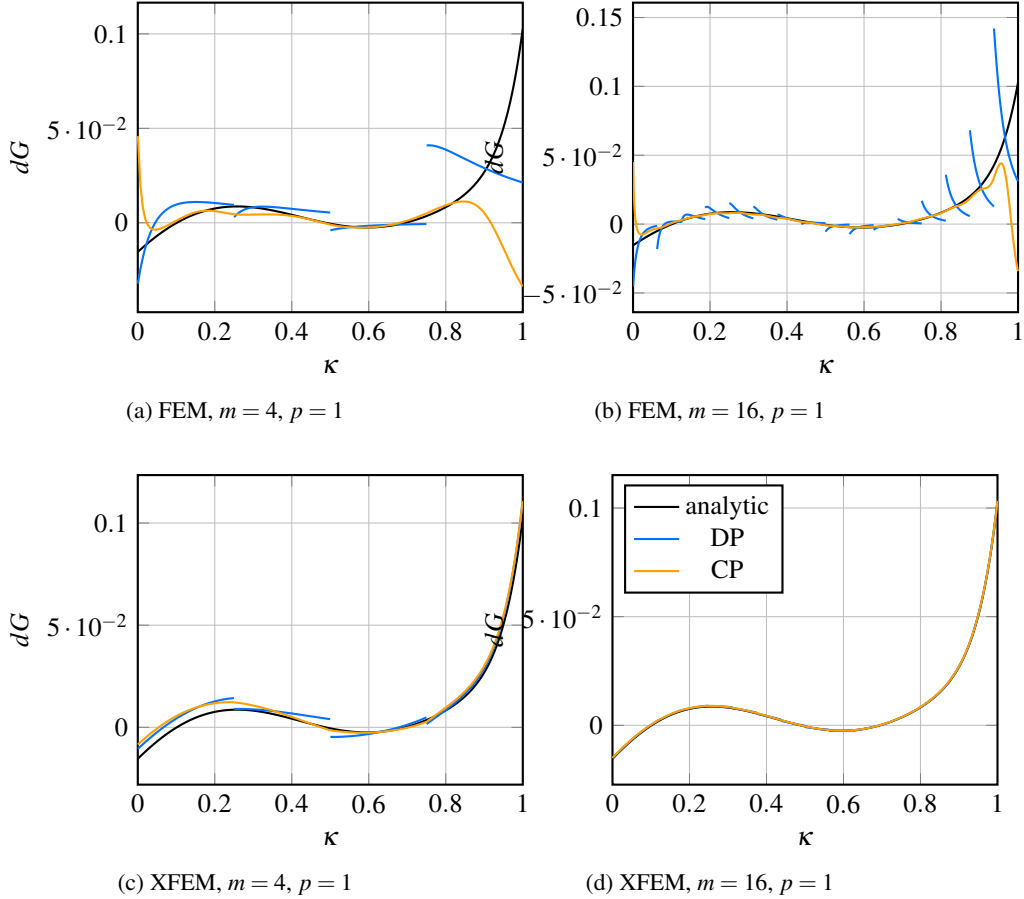
In this section, we analyze the differences of the discretize-then-optimize and optimize-then-discretize approaches for the shape derivative, \corrections{i.e., we compare the formulas for the numerical shape derivative obtained either with the \standardMethod \eqref{eq::derivativeTrackingShapeStandard} or the \enrichedMethod \eqref{eq::derivativeTrackingShapeEnriched} with the discretizations of the continuous shape derivative formula \eqref{eq::continousShapeDerivativeVolume}}.

\begin{figure}[htb]
\begin{subfigure}[t]{0.99\textwidth}
\centering
\begin{tikzpicture}
\begin{axis}[hide axis,
legend style={draw=white!15!black,legend cell align=left,legend columns=-1,xmin=10,xmax=50,ymin=0,ymax=0.4,}    ]
\addlegendimage{col1,mark=o}
\addlegendentry{DP $p=1$};
\addlegendimage{col2,mark=+}
\addlegendentry{DP  $p = 2$};
\addlegendimage{col3,mark=square}
\addlegendentry{DP $p = 3$};
\addlegendimage{col4,mark=asterisk}
\addlegendentry{CP $p=1$};
\addlegendimage{col5,mark=star}
\addlegendentry{CP $p = 2$};
\addlegendimage{col6,mark=triangle}
\addlegendentry{CP $p = 3$};
\end{axis}
\end{tikzpicture}
\end{subfigure}

\begin{subfigure}[t]{0.49\textwidth}
    \begin{tikzpicture}
    \begin{axis}[grid,xmode=log,ymode=log,
    xlabel={$m$}, ylabel={$L_2$ error},
        width=6.7cm,height=5.7cm]

    \addplot[dashed,domain=2:256,black] {0.1/x^0.45)}; 
    \addplot[col1,mark=o] table[col sep=comma,y index=1] {pic/shape_FEM_L2_3.txt};
    \addplot[col2,mark=+] table[col sep=comma,y index=2] {pic/shape_FEM_L2_3.txt};
    \addplot[col3,mark=square] table[col sep=comma,y index=3] {pic/shape_FEM_L2_3.txt};

    \addplot[col4,mark=asterisk] table[col sep=comma,y index=4] {pic/shape_FEM_L2_3.txt};
    \addplot[col5,mark=star] table[col sep=comma,y index=5] {pic/shape_FEM_L2_3.txt};
    \addplot[col6,mark=triangle] table[col sep=comma,y index=6] {pic/shape_FEM_L2_3.txt};

    \legend{{$O(h^{0.45})$}}
    \end{axis}
    \end{tikzpicture}
    \caption{\FEM}
    \label{fig:shapeFEML2}
\end{subfigure}
\begin{subfigure}[t]{0.49\textwidth}
    \begin{tikzpicture}
    \begin{axis}[grid,xmode=log,ymode=log,
    xlabel={$m$}, ylabel={$L_2$ error},
        width=6.7cm,height=5.7cm]

    \addplot[dashed,domain=2:256,black] {0.1/x^2)};     \addplot[domain=2:256,black] {0.005/x^(3))};

    \addplot[col1,mark=o] table[col sep=comma,y index=7] {pic/shape_FEM_L2_3.txt};
    \addplot[col2,mark=+] table[col sep=comma,y index=8] {pic/shape_FEM_L2_3.txt};
    \addplot[col3,mark=square] table[col sep=comma,y index=9] {pic/shape_FEM_L2_3.txt};

    \addplot[col4,mark=asterisk] table[col sep=comma,y index=10] {pic/shape_FEM_L2_3.txt};
    \addplot[col5,mark=star] table[col sep=comma,y index=11] {pic/shape_FEM_L2_3.txt};
    \addplot[col6,mark=triangle] table[col sep=comma,y index=12] {pic/shape_FEM_L2_3.txt};

    \legend{{$O(h^{2})$},{$O(h^{3})$}}

    \end{axis}
    \end{tikzpicture}
    \caption{\XFEM}
    \label{fig:shapeXFEML2}
\end{subfigure}

    \caption{Convergence of the shape derivative $L_2$-error.}
    \label{fig:shapeL2}
\end{figure}
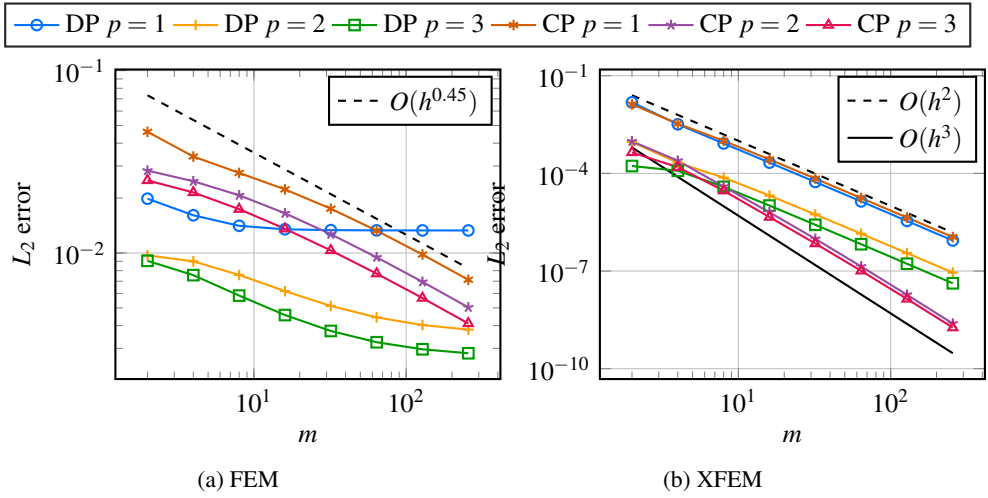

 In \Cref{fig:shapeFormulaA} the shape derivative obtained by \corrections{the \standardMethod} is plotted against the interface location. For the formula obtained by the discretize-then-optimize approach based on the discretized problem (DP) in \Cref{sec::discreteShapeStandard} we have already observed jumps in \Cref{fig:shapeFEMA}. In contrast to this the formula \eqref{eq::continousShapeDerivativeVolume} based on the continuous problem (CP) gives a continuous shape derivative. Furthermore, we observe that the CP approach matches the analytic solution within the interval $[0.2,0.8]$ well, however in the boundary regions large discrepancies can be observed. This behavior is also present if the mesh is refined, see \Cref{fig:shapeFormulaA2}. In contrast to this the results of the \corrections{\enrichedMethod (CP)} are smooth and also converge in the boundary regions, see \Cref{fig:shapeFormulaB,fig:shapeFormulaB2}. Moreover, while \corrections{the \enrichedMethod (DP)} yields a discontinuous result, the jumps disappear with finer mesh size. Finally, we look at the convergence rates of the shape derivative $L_2$ error defined by
 \begin{equation}
     e_{L_2} = \sqrt{\int_0^\ell (d \bar G(\kappa) - d_S \bar G(\kappa))^2 d\kappa}
 \end{equation}
 where $d \bar G(\kappa)$ is obtained by evaluating a finite element solution and $d_S \bar G(\kappa)$ is the analytic value of the shape derivative. The results for different discretizations and shape derivative formulas are given in \Cref{fig:shapeL2}. In \Cref{fig:shapeFEML2} the results for \corrections{the \standardMethod} are given. We observe no asymptotic convergence for when the DP formula \eqref{eq::derivativeTrackingShapeStandard} is used. In contrast to this for CP we observe a rate of approximately $0.45$. This low rate can be related to erroneous results in the boundary regions (see \Cref{fig:shapeFormulaA,fig:shapeFormulaA2}). The results for \corrections{the \enrichedMethod} are given in \Cref{fig:shapeXFEML2}. Here, we observe convergence for all variants. For $p=1$ the DP and the CP formulas give a rate of $2$. For $p=2,3$ the DP formula gives also a rate of $2$, whereas the CP formula gives a rate of $3$.

\subsection{Topological derivative}
\begin{figure}
    \begin{subfigure}[t]{0.47\textwidth}
		\centering
		\begin{tikzpicture}
        \begin{axis}[grid, xlabel={$\kappa$}, ylabel={$dG$},
            width=6.7cm,height=5.0cm, xmin=0,xmax=1,
            tick label style={/pgf/number format/fixed}]
        \addplot[only marks,col1,mark=o] table[col sep=comma,y index=1] {pic/top_FEM_1_n4_3.txt};
        \addplot[only marks,col2,mark=square] table[col sep=comma,y index=2] {pic/top_FEM_1_n4_3.txt};
        \addplot[only marks,col3,mark=+] table[col sep=comma,y index=1] {pic/top_XFEM_1_n4_3.txt};
        \addplot[black] table[col sep=comma,y index=1] {pic/top_analytic_solution_3.txt};
        \end{axis}
        \end{tikzpicture}
		\caption{$m=4$}
	\end{subfigure}
	\begin{subfigure}[t]{0.47\textwidth}
		\centering
		\begin{tikzpicture}
        \begin{axis}[grid, xlabel={$\kappa$}, ylabel={$dG$},
            width=6.7cm,height=5.0cm, xmin=0,xmax=1,
            tick label style={/pgf/number format/fixed}]
        \addplot[only marks,col1,mark=o] table[col sep=comma,y index=1] {pic/top_FEM_1_n8_3.txt};
        \addplot[only marks,col2,mark=square] table[col sep=comma,y index=2] {pic/top_FEM_1_n8_3.txt};
        \addplot[only marks,col3,mark=+] table[col sep=comma,y index=1] {pic/top_XFEM_1_n8_3.txt};
        \addplot[black] table[col sep=comma,y index=1] {pic/top_analytic_solution_3.txt};
        \end{axis}
        \end{tikzpicture}
		\caption{$m=8$}
	\end{subfigure}
	\vspace{0.5cm}

    \begin{subfigure}[t]{0.47\textwidth}
		\centering
		\begin{tikzpicture}
        \begin{axis}[grid, xlabel={$\kappa$}, ylabel={$dG$},
            width=6.7cm,height=5.0cm, xmin=0,xmax=1,
            tick label style={/pgf/number format/fixed}]
        \addplot[only marks,col1,mark=o] table[col sep=comma,y index=1] {pic/top_FEM_1_n16_3.txt};
        \addplot[only marks,col2,mark=square] table[col sep=comma,y index=2] {pic/top_FEM_1_n16_3.txt};
        \addplot[only marks,col3,mark=+] table[col sep=comma,y index=1] {pic/top_XFEM_1_n16_3.txt};
        \addplot[black] table[col sep=comma,y index=1] {pic/top_analytic_solution_3.txt};
        \end{axis}
        \end{tikzpicture}
		\caption{$m=16$}
	\end{subfigure}
	\begin{subfigure}[t]{0.47\textwidth}
		\centering
		\begin{tikzpicture}
        \begin{axis}[grid, xlabel={$\kappa$}, ylabel={$dG$},
            width=6.7cm,height=5.0cm, xmin=0,xmax=1,
            tick label style={/pgf/number format/fixed}]
        \addplot[only marks,col1,mark=o] table[col sep=comma,y index=1] {pic/top_FEM_1_n32_3.txt};
        \addplot[only marks,col2,mark=square] table[col sep=comma,y index=2] {pic/top_FEM_1_n32_3.txt};
        \addplot[only marks,col3,mark=+] table[col sep=comma,y index=1] {pic/top_XFEM_1_n32_3.txt};
        \addplot[black] table[col sep=comma,y index=1] {pic/top_analytic_solution_3.txt};
        \end{axis}
        \end{tikzpicture}
		\caption{$m=32$}
	\end{subfigure}

    \begin{subfigure}[t]{0.99\textwidth}
		\centering
		\begin{tikzpicture}
\begin{axis}[hide axis,
legend style={draw=white!15!black,legend cell align=left,legend columns=-1,xmin=10,xmax=50,ymin=0,ymax=0.4,}    ]
\addlegendimage{black}
\addlegendentry{analytic};
\addlegendimage{col1,mark=o}
\addlegendentry{\FEM};
\addlegendimage{col3,mark=+}
\addlegendentry{\XFEM};
\addlegendimage{col2,mark=square}
\addlegendentry{\textit{corrected} \standardMethod};
\end{axis}
\end{tikzpicture}
	\end{subfigure}
	\caption{Topological derivative for different discretizations}
	\label{fig:tracking_results_topDerivative}
\end{figure}
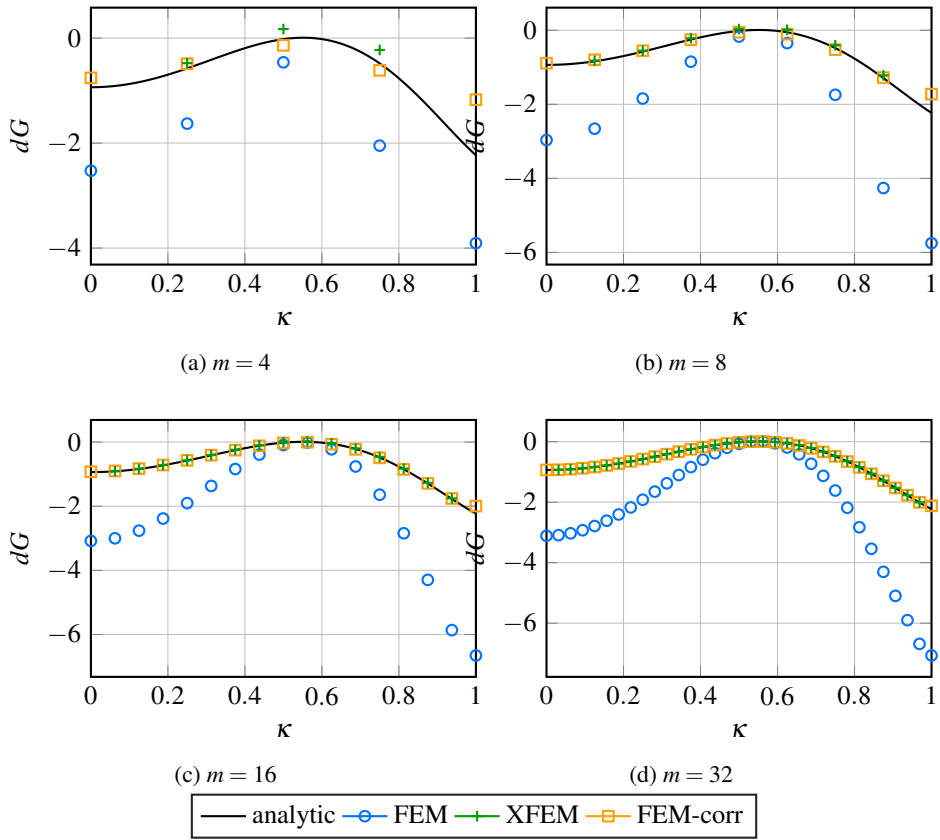
In \Cref{fig:tracking_results_topDerivative}, the results for the topological derivative for different numbers of elements are presented. We observe that, for the \standardMethod, the topological derivatives computed using formula \eqref{eq::topologicalDerivativeStandard} do not converge to the analytical continuous topological derivatives. However, when formula \eqref{eq::topologicalDerivativeStandard} is multiplied by the additional factor of $\frac{\lambda_1}{\lambda_2}$ (denoted as \textit{corrected} \standardMethod), the correct results are obtained. This factor is motivated by comparing the numerical sensitivity \eqref{eq::topologicalDerivativeStandard} with the analytical one \eqref{eq_ana_TD}. The discrepancy comes from the fact that the polarization matrix is not correctly accounted for in the discretized case if kinks of the solution across material interfaces are not resolved; see also the discussion in \cite[Rem. 6]{gangl2023unified}.

In contrast, for the \enrichedMethod, highly accurate topological derivatives are obtained directly using formula \eqref{eq::topological_derivativeC}, without the need for further modifications.

\subsection{Effect of material contrast}
The numerical results presented in this section were given for the relatively small contrast of a factor of 3 between $\lambda_1=0.6$ and $\lambda_2=0.2$. While for most realistic topology optimization problems, the contrast is much higher (typically chosen as $10^4$), we chose these values for a clearer presentation of the observed discontinuities and oscillations. In this subsection, we comment on the effects appearing in the case of material parameter values, $\lambda_2=0.2$ and $\lambda_1 = 10^6 \lambda_2$, and show that the main conclusions drawn above persist. For brevity, we focus on the sensitivities of the discretized problem (DP), illustrated in \Cref{fig:shapeFEM,fig:shapeXFEM,fig:shapeL2}.

\Cref{fig:contrastShapeFEM} shows the shape derivative in the case of the \standardMethod as a function of the interface position for different polynomial degrees $p$, as in \Cref{fig:shapeFEM}. We observe that the shape derivative is highly oscillatory for $p=1,2,3$ with large deviation from the analytic shape derivative. The shape derivative obtained by the \enrichedMethod in this high-contrast setting is depicted in \Cref{fig:contrastShapeXFEM}. Here, when compared with the low-contrast case \Cref{fig:shapeXFEM}, the oscillations for $p=2,3$ become more evident. However, as the number of mesh nodes $m$ increases, the shape derivative using the \enrichedMethod still converges whereas it does not in the case of the \standardMethod, see \Cref{fig:shapeL2highContrast}. Finally, we mention that, in the case of the topological derivative, in contrast to the results obtained by the \standardMethod, the results obtained by the \enrichedMethod perfectly match the analytic solution, see \Cref{fig:ContrastTop}. Again, the same holds true for the \textit{corrected} \standardMethod.

\begin{figure}[htb]
\begin{subfigure}[t]{0.49\textwidth}
    \begin{tikzpicture}
    \begin{axis}[grid, xlabel={$\kappa$}, ylabel={$dG$},
        width=6.7cm,height=5.7cm,
        xmin=0,xmax=1,
        restrict y to domain = -0.05:1.8,
        ymax = 0.1, ymin = -0.05,
        unbounded coords=jump,tick label style={/pgf/number format/fixed}]

    \addplot[black] table[col sep=comma,y index=2] {pic/shape_analytic_solution_1000000.txt};
    \addplot[col1] table[col sep=comma,y index=2] {pic/shape_FEM_1_n16_1000000.txt};
    \addplot[col2] table[col sep=comma,y index=2] {pic/shape_FEM_2_n16_1000000.txt};
    \addplot[col3] table[col sep=comma,y index=2] {pic/shape_FEM_3_n16_1000000.txt};
    \end{axis}
    \end{tikzpicture}
    \caption{Shape derivative for \standardMethod.}
    \label{fig:contrastShapeFEM}
\end{subfigure}
\begin{subfigure}[t]{0.49\textwidth}
\centering
    \begin{tikzpicture}
    \begin{axis}[grid, xlabel={$\kappa$}, ylabel={$dG$},
        width=6.7cm,height=5.7cm,xmin=0,xmax=1,
        ymax = 0.1, ymin = -0.05,
        restrict y to domain = -0.05:1.8, legend style={font=\scriptsize, at={(0.4,0.05)},anchor=south}, unbounded coords=jump,
        tick label style={/pgf/number format/fixed} ]

    \addplot[black] table[col sep=comma,y index=2] {pic/shape_analytic_solution_1000000.txt};
    \addplot[col1] table[col sep=comma,y index=2] {pic/shape_XFEM_1_n16_1000000.txt};
    \addplot[col2] table[col sep=comma,y index=2] {pic/shape_XFEM_2_n16_1000000.txt};
    \addplot[col3] table[col sep=comma,y index=2] {pic/shape_XFEM_3_n16_1000000.txt};
    \legend{{analytic},$p=1$,$p=2$,$p=3$}
    \end{axis}
    \end{tikzpicture}
    \caption{Shape derivative for \enrichedMethod.}
    \label{fig:contrastShapeXFEM}
\end{subfigure}
\vspace{0.5cm}

\begin{subfigure}[t]{0.49\textwidth}
\begin{tikzpicture}
    \begin{axis}[grid,xmode=log,ymode=log,
    xlabel={$m$}, ylabel={$L_2$ error},
        width=6.7cm,height=5.7cm, legend pos=south west, legend style={font=\footnotesize}]

    \addlegendimage{black,mark=o}
    \addlegendentry{ $p=1$};
    \addlegendimage{black,mark=+}
    \addlegendentry{ $p = 2$};
    \addlegendimage{black,mark=square}
    \addlegendentry{$p = 3$};
    \addlegendimage{black,dashed}
    \addlegendentry{$O(h^2)$};

    \addplot[col1,mark=o] table[col sep=comma,y index=1] {pic/shape_FEM_L2_1000000.txt};

    \addplot[col1,mark=+] table[col sep=comma,y index=2] {pic/shape_FEM_L2_1000000.txt};
    \addplot[col1,mark=square] table[col sep=comma,y index=3] {pic/shape_FEM_L2_1000000.txt};

    \addplot[dashed,domain=2:256,black] {0.08/x^2)};

    \addplot[col3,mark=o] table[col sep=comma,y index=7] {pic/shape_FEM_L2_1000000.txt};
    \addplot[col3,mark=+] table[col sep=comma,y index=8] {pic/shape_FEM_L2_1000000.txt};
    \addplot[col3,mark=square] table[col sep=comma,y index=9] {pic/shape_FEM_L2_1000000.txt};

    \end{axis}
    \end{tikzpicture}

    \caption{$L^2$ error of shape derivatives.}
    \label{fig:shapeL2highContrast}
\end{subfigure}
\begin{subfigure}[t]{0.47\textwidth}
		\centering
		\begin{tikzpicture}
        \begin{axis}[grid, xlabel={$\kappa$}, ylabel={$dG$},
            width=6.7cm,height=5.7cm, xmin=0,xmax=1,legend pos=south west, legend style={font=\scriptsize}]

        \addlegendimage{black}
        \addlegendentry{analytic};
        \addlegendimage{col1,mark=o}
        \addlegendentry{\FEM};
        \addlegendimage{col3,mark=+}
        \addlegendentry{\XFEM};
        \addlegendimage{col2,mark=square}
        \addlegendentry{\textit{corr.} \standardMethod};

        \addplot[only marks,col1,mark=o] table[col sep=comma,y index=1] {pic/top_FEM_1_n32_1000000.txt};
        \addplot[only marks,col2,mark=square] table[col sep=comma,y index=2] {pic/top_FEM_1_n32_1000000.txt};
        \addplot[only marks,col3,mark=+] table[col sep=comma,y index=1] {pic/top_XFEM_1_n32_1000000.txt};
        \addplot[black] table[col sep=comma,y index=1] {pic/top_analytic_solution_1000000.txt};
        \end{axis}
        \end{tikzpicture}
		\caption{Topological derivative}
        \label{fig:ContrastTop}
	\end{subfigure}
\begin{subfigure}[t]{0.45\textwidth}
		\centering
		\begin{tikzpicture}
\begin{axis}[hide axis,
legend style={draw=white!15!black,legend cell align=left,legend columns=1,xmin=10,xmax=50,ymin=0,ymax=0.4,}    ]

\end{axis}
\end{tikzpicture}
	\end{subfigure}

    \caption{\corrections{Numerical results for high contrast $\lambda_1/\lambda_2 = 10^6$: (a) Shape derivative for \standardMethod for $m=16$ and $p=1,2,3$, cf. Figure \ref{fig:shapeFEMB}; (b) Shape derivative for \enrichedMethod for $m=16$ and $p=1,2,3$, cf. Figure \ref{fig:shapeXFEMB}; (c) $L^2$ error for \standardMethod (blue) and \enrichedMethod (green), cf. Figure \ref{fig:shapeL2}; (d) Topological derivative for \standardMethod, \enrichedMethod and \textit{corrected} \standardMethod, cf. Figure \ref{fig:tracking_results_topDerivative}.}}
    \label{fig:highContrast}
\end{figure}

\section{Conclusion}
We have studied the influence of the discretization on the shape and topological derivative for a two material optimization problem with an elliptic PDE constraint in 1D. In particular we studied how the smoothness of the shape functions affects the shape derivative. To conduct this study we used splines as smooth shape functions. It turns out that the smoothness of the shape derivative is related to the smoothness of the discretization. In particular for linear finite elements ($p=1$) the shape derivative has jumps at the element interfaces, which might spoil the performance of gradient based optimization algorithms. For higher order and thus smoother splines the shape derivative is at least continuous.

Furthermore, we studied how the interface resolution by the discretization affects the shape and the topological derivative. If no interface position is incorporated in the discretization the shape derivative oscillates which might also spoil the performance of gradient based optimization algorithms. Furthermore, the topological derivative does not converge to the analytic solution. In contrast to this an enhanced discretization, which takes the interface position into account, yields converging shape and topological derivatives.

As a third aspect we compared the discretize-then-optimize and the optimize-then-discretize approaches. Linear finite elements and the optimize-then-discretize approach yield a continuous and less oscillatory shape derivative. However, this assertion is applicable exclusively to the "interior" of the domain. It does not extend to the boundary regions where erroneous shape derivative values are computed. For enriched discretizations we observe higher convergence rates for the optimize-then-discretize approach, giving shape derivatives closer to the shape derivative of the continuous problem.

To summarize, we have seen a superior behavior of the enhanced discretizations taking the interface position into account, which gives fast converging non-oscillatory sensitivities. Furthermore, smoother shape functions in the discretizations yield smoother shape derivatives, which might be also beneficial for gradient based optimization algorithms. Finally, we observed higher convergence rates for the optimize-then-discretize approach than for the discretize-then-optimize to the sensitivities of the continuous problem.

\corrections{
The present study has been restricted to a spatially one-dimensional setting. This allowed us to directly compute analytic expressions for the shape and topological derivative, it significantly alleviated the computation of numerical shape and topological derivatives and allowed to visualize them as functions of the interface location. In a $d$-dimensional setting ($d=2,3$), the interface will be a $(d-1)$-dimensional manifold and an illustrative presentation of the results will be much harder. Furthermore, the number of enrichment functions $n_E$ needed for accurately resolving the material interfaces will be much larger and will vary in the course of an optimization process. Nevertheless, we expect that the sensitivity analysis in this higher-dimensional setting can be carried out along the lines of the present work.

}

\paragraph{Acknowledgment}
The work of P.G. is partially supported by the joint DFG/FWF Collaborative Research Centre CREATOR (DFG Project-ID 492661287/TRR 361; FWF Project-DOI 10.55776/F90) at TU Darmstadt, TU Graz and JKU/RICAM Linz. For open access purposes, the author has applied a CC BY public copyright license to any author accepted manuscript version arising from this submission. Moreover, P.G. is partially supported by the State of Upper Austria.

\begin{appendix}
\renewcommand*{\thesection}{\Alph{section}}
	\section{Analytic expression of the objective function}\label{appandix::analyticG}
	The objective function for $f(x)=x$ and $\bar u(x)=x(\ell-x)$ can be written as the sum of six terms,
	
\begin{equation}\label{eq::G_kappa_explicit}
   G(\kappa) = G_1(\kappa) + G_2(\kappa) + G_3(\kappa)+ G_4(\kappa) + G_5(\kappa) + G_6(\kappa)
\end{equation}
 with
\begin{align*}
   G_1(\kappa) &= 
\frac{{a_{1}}^2\kappa ^7}{7}+\frac{2a_{1}a_{2}\kappa ^5}{5}+\frac{{a_{2}}^2\kappa ^3}{3},\\
  G_2(\kappa) &= 
\frac{\kappa ^3\left(6\kappa ^2-15\kappa +10\right)}{30},\\
  G_3(\kappa) &= 
\frac{a_{1}\kappa ^6}{3}-\frac{2a_{1}\kappa ^5}{5}+\frac{a_{2}\kappa ^4}{2}-\frac{2a_{2}\kappa ^3}{3},\\
  G_4(\kappa) &= 
\frac{{b_{2}}^2\left({{\ell}}^3-\kappa ^3\right)}{3}+\frac{{b_{1}}^2\left({{\ell}}^7-\kappa ^7\right)}{7}+{b_{3}}^2\left({\ell}-\kappa \right) \nonumber \\ 
 \quad& +b_{2}b_{3}\left({{\ell}}^2-\kappa ^2\right)+\frac{b_{1}b_{3}\left({{\ell}}^4-\kappa ^4\right)}{2}+\frac{2b_{1}b_{2}\left({{\ell}}^5-\kappa ^5\right)}{5},\\
  G_5(\kappa) &= 
\frac{{{\ell}}^5}{5}-\frac{{{\ell}}^4}{2}+\frac{{{\ell}}^3}{3}-\frac{\kappa ^5}{5}+\frac{\kappa ^4}{2}-\frac{\kappa ^3}{3},\\
  G_6(\kappa) &= 
\frac{2b_{3}\left({{\ell}}^3-\kappa ^3\right)}{3}-\frac{2b_{2}\left({{\ell}}^3-\kappa ^3\right)}{3}-b_{3}\left({{\ell}}^2-\kappa ^2\right) \nonumber\\ 
 \quad& +\frac{b_{2}\left({{\ell}}^4-\kappa ^4\right)}{2}-\frac{2b_{1}\left({{\ell}}^5-\kappa ^5\right)}{5}+\frac{b_{1}\left({{\ell}}^6-\kappa ^6\right)}{3},
\end{align*}
 with
\begin{align*}
a_{1} &= -\frac{1}{6\lambda _{1}},\\
a_{2} &= \frac{{{\ell}}^3\lambda _{1}-\kappa ^3\lambda _{1}+\kappa ^3\lambda _{2}}{6\lambda _{1}\left({\ell}\lambda _{1}-\kappa \lambda _{1}+\kappa \lambda _{2}\right)},\\
b_{1} &= -\frac{1}{6\lambda _{2}},\\
b_{2} &= \frac{{{\ell}}^3\lambda _{1}-\kappa ^3\lambda _{1}+\kappa ^3\lambda _{2}}{6\lambda _{2}\left({\ell}\lambda _{1}-\kappa \lambda _{1}+\kappa \lambda _{2}\right)},\\
b_{3} &= \frac{{\ell}\left(\kappa ^3\lambda _{1}-\kappa ^3\lambda _{2}-{{\ell}}^2\kappa \lambda _{1}+{{\ell}}^2\kappa \lambda _{2}\right)}{6\lambda _{2}\left({\ell}\lambda _{1}-\kappa \lambda _{1}+\kappa \lambda _{2}\right)}.
\end{align*}

\end{appendix}

\bibliographystyle{plainnat}
\bibliography{literature}

\end{document}